\documentclass[11pt,a4paper,reqno,twoside]{amsart}

\usepackage[T1]{fontenc}
\usepackage[utf8]{inputenc}
\usepackage[english]{babel}
\usepackage{amsmath,amssymb,amsthm,mathtools,mathrsfs}
\usepackage{enumitem}
\usepackage[colorlinks=true,linkcolor=blue,citecolor=blue,urlcolor=blue]{hyperref}
\usepackage[nameinlink,capitalise,noabbrev]{cleveref}

\allowdisplaybreaks
\numberwithin{equation}{section}
\usepackage{amsthm}
\usepackage{aliascnt}

\newtheorem{theorem}{Theorem}[section]

\newaliascnt{proposition}{theorem}

\aliascntresetthe{proposition}

\newaliascnt{lemma}{theorem}
\newtheorem{lemma}[lemma]{Lemma}
\aliascntresetthe{lemma}

\newaliascnt{corollary}{theorem}

\aliascntresetthe{corollary}

\newaliascnt{definition}{theorem}
\newtheorem{definition}[definition]{Definition}
\aliascntresetthe{definition}

\newaliascnt{remark}{theorem}
\newtheorem{remark}[remark]{Remark}
\aliascntresetthe{remark}

\usepackage{hyperref}
\usepackage[nameinlink,noabbrev,capitalise]{cleveref}

\crefname{theorem}{theorem}{theorems}
\Crefname{theorem}{Theorem}{Theorems}

\crefname{proposition}{proposition}{propositions}
\Crefname{proposition}{Proposition}{Propositions}

\crefname{lemma}{lemma}{lemmas}
\Crefname{lemma}{Lemma}{Lemmas}

\crefname{corollary}{corollary}{corollaries}
\Crefname{corollary}{Corollary}{Corollaries}

\crefname{definition}{definition}{definitions}
\Crefname{definition}{Definition}{Definitions}

\crefname{remark}{remark}{remarks}
\Crefname{remark}{Remark}{Remarks}

\newcommand{\R}{\mathbb R}
\newcommand{\norm}[1]{\left\lVert#1\right\rVert}
\newcommand{\dd}{\,\mathrm d}
\newcommand{\eps}{\varepsilon}
\newcommand{\K}{\mathcal K}
\newcommand{\E}{\mathcal E}
\newcommand{\N}{\mathcal N}

\title[Exact Mass Threshold for anti-de Sitter Equations]{An Exact Mass Threshold for Global Existence and Blow-up\\
in Semilinear Anti-de Sitter Equations}

\author{Wanderley Nunes do Nascimento$^{1}$}
\address{$^{1}$ Department of Pure and Applied Mathematics, Federal University of Rio Grande do Sul, Porto Alegre, Brazil}
\author{Michael Reissig$^2$}
\address{$^{2}$Institute of Applied Analysis, Faculty of Mathematics and
Computer Science, TU Bergakademie Freiberg, Pr\"uferstra\ss e 9,
09596 Freiberg, Germany.}

\begin{document}

\begin{abstract}
We study the influence of the mass parameter on global existence and
finite-time blow-up for the semilinear anti-de Sitter Klein--Gordon
equation
\[
 \phi_{tt}-e^{2t}\Delta\phi+n\phi_t+m^2\phi=|\phi|^p
 \qquad\text{in }[0,\infty)\times\R^d.
\]
After removing the damping term, we establish $L^2$ energy estimates for
the resulting Klein--Gordon equation with constant coefficient in the potential. They imply
global existence for arbitrary sufficiently small energy data, without
spatial localization, whenever
\(
 m^2>-\frac{d(2n+d)}{4}
 \quad\text{and}\quad p>1.
\)
In the remaining mass parameter interval above the structural threshold, an exact
change to conformal time reduces the problem to a scale-invariant
effective damped wave equation. Exponentially weighted energy estimates
then extend global existence to
\(
  -d(n+d)<m^2\le-\frac{d(2n+d)}4
 \quad\text{for every } \quad p>1.
\)
In the complementary regime $m^2\le-d(n+d)$, we prove finite-time blow-up for every $p>1$ and for weak solutions with finite speed of propagation arising from nonnegative, nontrivial initial data. We further obtain lifespan estimates that are polynomial at the threshold and logarithmic below it. Thus, $m^2=-d(n+d)$ is the exact structural mass parameter threshold separating the global-existence and blow-up regimes.
\end{abstract}
\maketitle
\vspace{-2em}

\section{Introduction}

The propagation of waves on cosmological backgrounds is a natural meeting
point between geometric analysis, mathematical physics, and the qualitative
theory of nonlinear hyperbolic equations. Two of the basic constant-curvature
models are de Sitter and anti-de Sitter spacetimes, distinguished by the sign
of the cosmological constant. De Sitter spacetime has positive cosmological
constant and exhibits accelerated expansion, whereas anti-de Sitter
spacetime has negative cosmological constant and, in the coordinates relevant
to the present work, leads to a contracting geometry. For spatially
homogeneous metrics, these geometric effects appear in the reduced wave
equation through a time-dependent propagation speed and a first-order term
generated by the variation of the spatial volume. Thus, the expanding de
Sitter model produces an exponentially decreasing propagation speed, while
its contracting anti-de Sitter counterpart produces the exponentially
increasing speed considered below. Fundamental solutions and integral
representations of solutions to wave and Klein--Gordon equations on de Sitter,
anti-de Sitter, and more general Robertson--Walker backgrounds were developed
in \cite{Yagdjian2010,Yagdjian2015,YagdjianGalstian2008,YagdjianGalstian2009dS,
YagdjianGalstian2009AdS}; see also
\cite{Baskin2013,HintzVasy2015} for dispersive and microlocal approaches on
asymptotically de Sitter spaces.

The Klein--Gordon equation is the canonical model for a massive scalar field
on a Lorentzian manifold. On curved backgrounds, Klein--Gordon-type equations
also arise in the study of scalar perturbations, symmetry breaking,
gravitational and electromagnetic perturbations, Kaluza--Klein reductions,
and models from high-energy physics. In the global anti-de Sitter geometry,
an additional subtlety is caused by the timelike character of the conformal
boundary, so that the global dynamics generally depends on the boundary
condition imposed at conformal infinity; see
\cite{Holzegel2012,IshibashiWald2004,Vasy2012,Warnick2013}.
The model considered in the present paper, however, is posed as a Cauchy
problem on $\mathbb R^d$ and does not involve a boundary condition at
conformal infinity.

The influence of a negative mass squared also depends strongly on the
underlying geometry. For Klein--Gordon equations on an expanding de Sitter background, an
imaginary mass gives rise to a linearly unstable regime in which the linear
dynamics admits exponentially growing modes, with important consequences
for the corresponding nonlinear problem; see, for instance,
\cite{Yagdjian2009,YagdjianGalstian2009dS}. In the anti-de Sitter setting,
by contrast, a negative mass squared does not by itself imply linear
instability. The classical Breitenlohner--Freedman theory shows that a
suitable renormalized energy remains positive above a distinguished negative
mass bound
\cite{BreitenlohnerFreedman1982a,BreitenlohnerFreedman1982b}.
This principle underlies a substantial part of the modern theory of massive
waves on asymptotically anti-de Sitter backgrounds; related well-posedness,
boundedness, and coupled Einstein--Klein--Gordon results can be found in
\cite{HolzegelSmulevici2012,HolzegelWarnick2014,Warnick2013}.
It also makes the mass parameter a genuine dynamical quantity: rather than
merely shifting the linear spectrum, it may determine whether the nonlinear
evolution is globally stable or develops singular behavior in finite time.

There is by now an extensive literature on linear and semilinear
Klein--Gordon equations in de Sitter spacetime. Explicit kernels and integral
transform methods were established in
\cite{Yagdjian2010,Yagdjian2013,Yagdjian2015,YagdjianGalstian2009dS}.
Global existence and blow-up for semilinear scalar fields were studied by
Yagdjian \cite{Yagdjian2009,Yagdjian2012,Yagdjian2019}, while energy
solutions and nonlinear models with power or quartic interactions were
treated in \cite{Nakamura2014,Nakamura2015,Nakamura2021}. Small-data global
existence based on regularity and dispersive estimates was obtained by
Ebert and Reissig \cite{EbertReissig2018}, and semilinear equations on
asymptotically de Sitter, Kerr--de Sitter, and related spacetimes were
investigated by Hintz and Vasy \cite{HintzVasy2015}. We also mention the
blow-up and lifespan analysis in de Sitter and Einstein--de Sitter models
developed in \cite{GalstianYagdjian2020,TsutayaWakasugi2020,
TsutayaWakasugi2022}.

The corresponding semilinear theory for the contracting anti-de Sitter model
is considerably less complete. At the linear level, Galstian
\cite{Galstian2003} obtained $L^p$--$L^q$ estimates for equations with
exponentially increasing propagation speed, Yagdjian and Galstian
\cite{YagdjianGalstian2009AdS} derived representation formulas adapted to
the anti-de Sitter geometry, and Bachelot \cite{Bachelot2011} established
finite-energy, decay, and Strichartz estimates in the Poincar\'e chart.
For the nonlinear problem, Palmieri and Takamura
\cite{PalmieriTakamura2022} proved blow-up results for de Sitter and
anti-de Sitter models by combining lower bounds for spatial averages with an
iteration procedure. They subsequently treated the critical anti-de Sitter
case in \cite{PalmieriTakamura2022critical}. These works demonstrate the
strong influence of the exponentially growing propagation speed, but they do
not identify a mass value separating global existence from blow-up for the
pure power nonlinearity. Determining such a threshold is the main purpose of
the present paper.

Motivated by this framework, we consider the Cauchy problem
\begin{equation}\label{eq:phi}
 \begin{cases}
 \phi_{tt}-e^{2t}\Delta\phi+n\phi_t+m^2\phi=|\phi|^p,
 &(t,x)\in[0,\infty)\times\R^d,\\
 \phi(0,x)=f(x),\qquad \phi_t(0,x)=g(x),&x\in\R^d,
 \end{cases}
\end{equation}
where $d\ge1$ is the spatial dimension, $n>0$ is the coefficient of the
dissipative term, $m^2\in\R$, and $p>1$. In the geometric Klein--Gordon
equation the damping coefficient is tied to the spatial dimension. Here we
keep $d$ and $n$ independent in order to isolate the competing effects of
geometry and dissipation: $d$ governs the exponential growth of the volume
of the domain of dependence, whereas $n$ determines the damping and the
location of the characteristic roots.

Our objective is to determine how the mass interacts with the dissipation and
with the contracting geometry, and in particular to identify a single mass
threshold separating small-data global existence from finite-time blow-up.
For $m^2<0$, the relevant linear growth exponent is
\begin{equation}\label{eq:rho}
 \rho:=\sqrt{\frac{n^2}{4}-m^2}-\frac n2.
\end{equation}
The spatially homogeneous solution of the linear equation grows like
$e^{\rho t}$.
On the other hand, finite propagation confines a solution issued from
compactly supported data to a ball of radius comparable with $e^t$, whose
volume is therefore of order $e^{dt}$. The balance between these two effects
occurs at $\rho=d$, or equivalently at
\begin{equation}\label{eq:threshold}
 m^2=-d(n+d).
\end{equation}

Our main result shows that \eqref{eq:threshold} is the exact structural
threshold. The global-existence region is covered by two complementary
arguments. In the upper mass parameter range
\[
 m^2>-\frac{d(2n+d)}4,
\]
we remove the damping and perform a localized Fourier analysis of the
resulting Klein--Gordon propagator. The shrinking low-frequency region has
volume of order $e^{-dt}$ and yields an additional $e^{-\frac{d}{2}t}$ gain in
$L^2$. This compensates for the growth of the spatially homogeneous component and
yields global energy solutions for arbitrary sufficiently small data, without
any spatial localization assumption. In the intermediate range
\[
 -d(n+d)<m^2\le-\frac{d(2n+d)}4,
\]
the preceding non-weighted argument no longer closes. We instead introduce
the conformal time $\tau=e^t-1$ and conjugate the unstable linear mode. The
equation is thereby transformed exactly into a wave equation with
scale-invariant effective damping and a time-dependent nonlinearity.
Exponentially weighted energy estimates, in the spirit of
\cite{IkehataTanizawa,Palmieri2018}, then yield global existence for sufficiently small
Gaussian-localized data. Together, the two methods cover the full region
$m^2>-d(n+d)$ for every power in the energy-subcritical range.

The complementary regime is governed by the opposite side of the same
balance. If
\[
 m^2\le-d(n+d),
\]
we integrate the equation in space and remove the mass term through the
exponential weight associated with the unstable characteristic root. Finite
propagation and H\"older's inequality then produce a weighted nonlinear
differential inequality for the spatial average, whose decisive coefficient
is $e^{(p-1)(\rho-d)t}$. Since $\rho\ge d$ precisely in this mass parameter range, an
ODE argument proves finite-time blow-up for every $p>1$ for weak solutions
with finite speed of propagation arising from nonnegative, nontrivial,
compactly supported data. The same reduction gives a polynomial lifespan
upper bound at the threshold and a logarithmic one below it.

The novelty is therefore not merely a collection of global-existence and
blow-up statements obtained in disjoint mass parameter regimes. The three
arguments reveal complementary manifestations of the same mechanism and
meet at one explicit mass value. At $m^2=-d(n+d)$, the exponential growth of the spatially homogeneous
component is exactly balanced by the geometric growth of the domain of
dependence. Above this value, dispersion or weighted dissipation controls
the power nonlinearity;
at and below it, the corresponding weighted spatial average necessarily
blows up. In this structural sense, the threshold is independent of $p$:
within the natural energy range on the existence side, the mass selects the
qualitative behavior, while the power determines the nonlinear time scale.

The paper is organized as follows. Section~2 states the global-existence and
blow-up results. Section~3 develops suitable estimates for Sobolev solutions to parameter-dependent Cauchy problems, the conformal-time
reduction and weighted energy method, and the spatial-average argument used
to prove the three main theorems.

\section{Main results}
We now state the global-existence and blow-up results that establish the
mass threshold described in the introduction. We begin by introducing the
parameters and the range of powers appearing in their formulations. Let
\begin{equation}\label{eq:q}
 q_m:=m^2-\frac{n^2}{4}.
\end{equation}
For every value of $q_m$, we set
\begin{equation}\label{parameternu}
 \nu:=\sqrt{|q_m|}.
\end{equation}
We also introduce
\begin{equation}\label{def:beta}
 \beta_0:=
 \begin{cases}
  \displaystyle
  \max\left\{\frac12,\nu\right\},
  & m^2<\dfrac{n^2}{4},\\[3mm]
  \dfrac12,
  & m^2\ge\dfrac{n^2}{4}.
 \end{cases}
\end{equation}
The parameter $\beta_0$ describes the maximal growth rate of the linear
propagator.

When $m^2<\frac{n^2}{4}$, we further introduce the
exponent
\begin{equation}\label{def:rho}
 \rho:=\nu-\frac n2.
\end{equation}
In particular, if $m^2<0$, then $\rho>0$ is the effective growth
exponent, and $e^{\rho t}$ describes the
growth of the unstable zero-frequency mode in the original variable
$\phi$.

%For the intermediate mass %regime, we also set
%\begin{equation}\label{def:kappa}
% \mu:=n+1+2\rho,
% \qquad
%\end{equation}
%These parameters arise naturally %from the exact reduction to a
%scale-invariant effectively %damped wave equation, while %$\Psi_0$
%determines the Gaussian weight %used in the second global-%existence
%theorem.

The global-existence results are proved in the energy-subcritical range
\begin{equation}\label{eq:p-range}
 \begin{cases}
  1<p<\infty, & d=1,2,\\[1mm]
  \displaystyle 1<p\le\frac{d}{d-2}, & d\ge3.
 \end{cases}
\end{equation}
The blow-up result, in contrast, holds for every $p>1$.
\begin{theorem}
\label{thm:global}
Let $d\ge1$, $n>0$, and
\begin{align}\label{RestmTheo1}
 m^2>-\frac{d(2n+d)}4.
\end{align}
Assume that $p>1$ satisfies \eqref{eq:p-range}.  There exists
$\varepsilon_0>0$ such that, if
\[
 (f,g)\in H^1(\R^d)\times L^2(\R^d),\,\,
\norm{f}_{H^1}+\norm{g}_{L^2}\le\varepsilon\le\varepsilon_0,
\]
then \eqref{eq:phi} has a unique global energy solution
\[
 \phi\in C\big([0,\infty),H^1(\R^d)\big)
 \cap C^1\big([0,\infty),L^2(\R^d)\big).
\]
Moreover, for every
\(
 \beta\in\left(\beta_0,\frac{n+d}{2}\right),
\)
where $\beta_0$ is defined in \eqref{def:beta}, the solution $\phi$
satisfies the energy estimate
\begin{equation}\label{eq:phi-energy-estimate}
 \norm{\phi(t)}_2+\norm{\phi_t(t)}_2
 +e^t\norm{\nabla\phi(t)}_2
 \le C_\beta\varepsilon
 e^{(\beta-\frac n2)t},
 \qquad t\ge0,
\end{equation}
where $C_\beta>0$ is independent of $t$ and $\varepsilon$.
\end{theorem}
In the following result $\Psi_0=\Psi_0(x)=\mu \frac{|x|^2}{2}$ with a suitable chosen $\mu>0$ in Section \ref{Section3.2}.
\begin{theorem}
\label{thm:global-weighted}
Let $d\ge1$ and $n \geq 1$, and let $p>1$ satisfy \eqref{eq:p-range}.
Assume that
\begin{equation}\label{rest:thm2}
 -d(n+d)<m^2\le -\frac{d(2n+d)}4,
\end{equation}
where the upper inequality is required to be strict when \(n=1\). There exists $\varepsilon_0>0$ such that, if
\[ \left(e^{\Psi_0}f,
 e^{\Psi_0}(g-\rho f) \right)\in  H^1(\R^d)\times L^2(\R^d) \,\,\,
\mbox{and}\,\,\,
\norm{e^{\Psi_0}f}_{H^1}
 +\norm{e^{\Psi_0}(g-\rho f)}_2
 \le\varepsilon\le\varepsilon_0,
\] where $\rho$ is defined in \eqref{def:rho},
then \eqref{eq:phi} has a unique global energy solution
\[
 \phi\in C\big([0,\infty),H^1(\R^d)\big)
 \cap C^1\big([0,\infty),L^2(\R^d)\big).
\]
Moreover, the solution $\phi$
satisfies the energy estimate
\begin{align}
 \norm{\phi(t)}_2 +\norm{\phi_t(t)}_2+e^t\norm{\nabla\phi(t)}_2
\le C\varepsilon e^{(\rho-\frac{d}{2})t}.
 \label{eq:weighted-phi-energy}
\end{align}
For every $2\le q\le \frac{2d}{(d-2)_+}\footnote{Here,
\((d-2)_+:=\max\{d-2,0\}\), and the quotient
\(2d/(d-2)_+\) is understood as \(+\infty\) when \(d=1\).
For \(d=2\), the endpoint \(q=\infty\) is excluded.}$ it holds
\begin{equation}\label{eq:weighted-phi-Lq}
 \norm{\phi(t)}_q
 \le C_q\varepsilon e^{(\rho-d(1-\frac{1}{q}))t}.
\end{equation}

\end{theorem}

\begin{remark}
The global existence conclusion remains also valid in the excluded
endpoint case \(n=1\) and \(m^2=-\frac{d(d+2)}{4}\). Indeed, the fixed-point argument can be carried out in the
corresponding logarithmically modified solution space. In this case, the
right-hand sides of \eqref{eq:weighted-phi-energy} and
\eqref{eq:weighted-phi-Lq} acquire the additional factor
\(
 1+\sqrt{t}.
\)
\end{remark}

% \begin{remark}
% Hypothesis \eqref{rest:thm2} implies that the parameter $\rho$ defined in \eqref{def:rho} satisfies
% \[
%  \frac d2\le\rho<d.
% \]
% \end{remark}

% \begin{remark}\label{rem:weighted-data}
% The Gaussian hypothesis is stronger than
% $(H^1\cap L^1)\times(L^2\cap L^1)$, but it is propagated by the
% scale-invariant dissipative flow and controls the $L^p$ norm required when
% $1<p<2$.  Thus the whole interval
% $-d(n+d)<m^2\le-d(2n+d)/4$ is covered.  The restriction
% \eqref{eq:p-range} is retained so that the
% nonlinearity is locally Lipschitz at the energy level.
% \end{remark}

Before stating the blow-up result, we introduce the notion of a weak
solution with finite speed of propagation that will be used throughout
the blow-up analysis. By $B_{R}$ we denote the ball around the origin with radius $R$ in $\R^d$.
\begin{definition}\label{def:weak-solution}
Let $T>0$, and assume that
\(
\operatorname{supp}f\cup\operatorname{supp}g\subset B_R.
\)
We say that
\(
 \phi\in L^p_{\mathrm{loc}}([0,T)\times\mathbb R^d)
\)
is a weak solution of \eqref{eq:phi} on $[0,T)$, with initial data
$(\varepsilon f,\varepsilon g)$, if the following two conditions are
satisfied:
\begin{enumerate}[label=\textnormal{(\roman*)}]
\item For almost every $t \in [0,T),$
\begin{align}\label{eq:weak-support}
 \operatorname{supp}\phi(t,\cdot)
 \subset B_{R+e^t-1};
\end{align}
\item for every
$\psi\in C_0^\infty([0,T)\times\mathbb R^d)$,
\begin{align}
 &\int_0^T\int_{\mathbb R^d}\phi
 \left(
   \psi_{tt}-e^{2t}\Delta\psi-n\psi_t+m^2\psi
 \right)\,dx\,dt
 \notag\\
 &\qquad=
 \int_0^T\int_{\mathbb R^d}|\phi|^p\psi\,dx\,dt
 +\varepsilon\int_{\mathbb R^d}(g+nf)\psi(0,x)\,dx
 -\varepsilon\int_{\mathbb R^d}f(x)\psi_t(0,x)\,dx.
 \label{eq:weak-formulation}
\end{align}
\end{enumerate}
The lifespan $T_\varepsilon$ is the supremum of all $T>0$ for which
such a weak solution exists.
\end{definition}

\begin{theorem}\label{thm:blowup}
Let $d\ge1$, $n\ge1$, $p>1$, and
\begin{align} \label{restonmublowup}
     m^2\le-d(n+d).
\end{align}
Assume that $f,g\in L^1(\R^d)$ are nonnegative, compactly supported, and not
both identically zero, with
\(
\operatorname{supp}f\cup\operatorname{supp}g\subset B_R.
\)
Then there exist $\varepsilon_0,C>0$, independent of $\varepsilon$, such that
every weak solution in the sense of \Cref{def:weak-solution} 
satisfies,
for $0<\varepsilon\le\varepsilon_0$,
\begin{equation}\label{eq:lifespan}
 T\le
 \begin{cases}
  C\varepsilon^{-(p-1)},&m^2=-d(n+d),\\[2mm]
  \displaystyle\frac{1}{\rho-d}\log(\varepsilon^{-1})+C,
       &m^2<-d(n+d),
 \end{cases}
\end{equation}
where $\rho$ is defined in \eqref{eq:rho}.
Consequently, the same estimate holds for $T_\varepsilon$, and in particular
$T_\varepsilon<\infty$.
\end{theorem}

\begin{remark}\label{rem:finite-propagation}
The support condition \eqref{eq:weak-support} is satisfied by energy
solutions constructed as mild solutions through the Duhamel formula.
Indeed, after introducing the conformal time variable
\(
 \tau=e^t-1,
\)
the principal part of \eqref{eq:phi} becomes
$\partial_\tau^2-\Delta$, while the damping, mass, and nonlinear terms
are local and of lower order. Hence, the standard domain-of-dependence
property gives
\[
 \operatorname{supp}\phi(t,\cdot)
 \subset B_{R+e^t-1}.
\]
Moreover, this property is preserved at every step of the fixed-point
iteration: the homogeneous linear solution has the stated support, and
each successive Duhamel term preserves the corresponding domain of
dependence.
In particular, when the initial data are compactly supported, the
energy solutions constructed in
\Cref{thm:global,thm:global-weighted} are mild solutions and belong to
the class introduced in \Cref{def:weak-solution}. Therefore,
\Cref{thm:blowup} constitutes the natural blow-up counterpart of the
global-existence results in
\Cref{thm:global,thm:global-weighted}.
\end{remark}

\section{Proofs of the Theorems}
\subsection{Proof of \Cref{thm:global}} \label{Section3.1}
% In order to prove \Cref{thm:global} we shall derive energy estimates for the correspondent linear equation. In order to get those estimates, we
% transform the Cauchy problem \eqref{eq:phi} we set
% \begin{equation}\label{eq:transform}
%  \phi(t,x)=e^{-nt/2}u(t,x).
% \end{equation}

Eliminating either the damping term or the mass term by a suitable change
of the unknown, and then applying Fourier analysis to the transformed
equation, is a standard strategy in the study of wave equations with
time-dependent lower-order coefficients; see, for instance,
\cite{NascimentoPalmieriReissig2017,Palmieri2018}.
In the present proof, we remove the damping term and derive suitable
estimates for solutions to the resulting homogeneous Klein--Gordon equation with
exponentially increasing propagation speed and constant coefficient in the potential. Our analysis is inspired by the approaches
developed in
\cite{BuiReissig2014,EbertNascimento2018,Galstian2003}. The Bessel-function
arguments used below follow the classical framework described in
\cite{WatsonBessel}.

We begin by introducing the transformation
\begin{equation}\label{eq:transform}
 \phi(t,x)=e^{-\frac n2t}u(t,x).
\end{equation}
A direct computation shows that the damping term is eliminated and
replaced by the potential $q_m u$ with
$q_m=m^2-\frac{n^2}{4}$
was defined in \eqref{eq:q}. Consequently, \eqref{eq:phi} is
equivalent to
\begin{equation}\label{eq:u}
 \begin{cases}
  u_{tt}-e^{2t}\Delta u+q_m u
  =e^{-\frac n2(p-1)t}|u|^p,
  &(t,x)\in[0,\infty)\times\R^d,\\
  u(0,x)=u_0(x):=f(x),&
  x\in\R^d,\\
  u_t(0,x)=u_1(x):=g(x)+\dfrac n2f(x),&
  x\in\R^d.
 \end{cases}
\end{equation}

We first study the linear equation associated with \eqref{eq:u}. Since its
principal coefficient \(e^{2t}\) depends explicitly on time, the
corresponding evolution is non-autonomous. Therefore, its propagators depend
on both the initial time \(s\) and the observation time \(t\). This
two-parameter formulation is also the natural one for Duhamel's principle,
because the nonlinear forcing generated at time \(s\) must be propagated
from \(s\) to \(t\).

More precisely, for \(t\ge s\ge0\), we denote by
\[
 K_0(t,s)
 \qquad\text{and}\qquad
 K_1(t,s)
\]
the propagators associated with the linear Cauchy problem
\begin{equation}\label{eq:linear}
 \begin{cases}
  v_{tt}-e^{2t}\Delta v+q_m v=0,\\
  v(s,x)=\varphi(x),\qquad
  v_t(s,x)=\psi(x).
 \end{cases}
\end{equation}
Thus, the corresponding solution is represented as
\[
 v(t)=K_0(t,s)\varphi+K_1(t,s)\psi.
\]
In particular, \(K_0(t,s)\) propagates the initial displacement, whereas
\(K_1(t,s)\) propagates the initial velocity.

For later use, we introduce the energy quantity
\begin{equation}\label{energe=1}
 \mathcal E_t(v)
 :=
 \norm{v(t)}_2+\norm{v_t(t)}_2
 +e^t\norm{\nabla v(t)}_2.
\end{equation}

The Fourier multipliers associated with \(K_0(t,s)\) and \(K_1(t,s)\)
will be expressed in terms of Bessel functions. We shall repeatedly use
their standard identities and asymptotic expansions, referring to
\cite{DLMF} for the specific formulas employed below.

% \begin{lemma}\label{lem:pure-L2}
% Set $\beta$ as defined in \eqref{def:beta}.
% For $t\ge s\ge0$,
% \begin{align}
%  \mathcal E_t(K_1(t,s)h)&\le C\Omega(t-s)e^{\beta(t-s)}\norm h_2,
%  \label{eq:pure-K1}\\
%  \mathcal E_t(K_0(t,s)h)&\le C\Omega(t-s)e^{\beta(t-s)}
%  (\norm h_2+e^s\norm{\nabla h}_2),\label{eq:pure-K0}
% \end{align}
% where $\Omega(\tau)=1$ if $q_m\neq0$, whereas
% $\Omega(\tau)\le C(1+\tau)$ in the borderline case $q_m=0$.
% Since $1+\tau\le C_\delta e^{\delta\tau}$ for every $\delta>0$,
% this factor can be absorbed into an arbitrarily small enlargement of
% $\beta$. Consequently,
% \begin{equation}\label{eq:pure-duhamel}
%  \sup_{0\le t\le T}e^{-\beta t}\mathcal E_t
%  \left(\int_0^tK_1(t,s)F(s)\,ds\right)
%  \le C\int_0^T e^{-\beta s}\norm{F(s)}_2\,ds.
% \end{equation}
% \end{lemma}

\begin{lemma}\label{lem:pure-L2}
Let $\beta_0$ and $\mathcal E_t$ be defined as in \eqref{def:beta} and \eqref{energe=1}. Then, for every
$t\geq s\geq0$, the propagators $K_0$ and $K_1$ of the Cauchy problem \eqref{eq:linear} satisfy
\begin{align}\mathcal E_t\bigl(K_0(t,s)h\bigr)
 &\leq
 C\Omega(t-s)e^{\beta_0(t-s)}
 \big(
  \norm{h}_{L^2}
  +e^s\norm{\nabla h}_{L^2}
 \big),
 \label{eq:pure-K0}\\
 \mathcal E_t\bigl(K_1(t,s)h\bigr)
 &\leq
 C\Omega(t-s)e^{\beta_0(t-s)}
 \norm{h}_{L^2},
 \label{eq:pure-K1}
 \end{align}
where
\[
 \Omega(\tau)=
 \begin{cases}
  1, & q_m\neq0,\\
  1+\tau, & q_m=0.
 \end{cases}
\]
Consequently, due to Duhamel's principle the solution of the Cauchy problem ($0\leq s \leq t\leq T$)
\begin{equation}\label{eq:linearinhomogeneous}
 v_{tt}-e^{2t}\Delta v+q_m v=F(t),
 \qquad
 v(s)=0,\quad v_t(s)=0
\end{equation}
satisfies the following estimates:\\
in the case $q_m\neq0$:
\begin{equation}\label{eq:pure-duhamel}
 \sup_{0\leq t\leq T}
 e^{-\beta_0 t}\mathcal E_t\Big(\int_0^t K_1(t,s)F(s)\,\dd s\Big)
 \leq
 C\int_0^T e^{-\beta_0 s}\norm{F(s)}_{L^2}\,\dd s.
\end{equation}
In the borderline case $q_m=0$, for every $\delta>0$,
\begin{equation}\label{eq:pure-duhamel-borderline}
 \sup_{0\leq t\leq T}
 e^{-(\beta_0+\delta)t}\mathcal E_t\Big(\int_0^t K_1(t,s)F(s)\,\dd s\Big)
 \leq
 C_\delta
 \int_0^T
 e^{-(\beta_0+\delta)s}\norm{F(s)}_{L^2}\,\dd s.
\end{equation}
\end{lemma}

\begin{proof}
 Throughout the proof, we use standard formulas for
Bessel functions collected in \cite[Chapter10]{DLMF}.
Taking the Fourier transform with respect to $x$ and writing
$r=|\xi|$, the equation of \eqref{eq:linear} becomes
\begin{equation}\label{eq:pure-fourier-linear}
 \widehat v_{tt}+(e^{2t}r^2+q_m)\widehat v=0.
\end{equation}
Introduce the new variable
\[
 \tau=e^tr
\]
and write $V(\tau)=\widehat v(t,\xi)$,
equation \eqref{eq:pure-fourier-linear} is transformed into
\begin{equation}\label{eq:pure-bessel-equation}
 \tau^2V''+\tau V'+(\tau^2+q_m)V=0.
\end{equation}
We first consider the case $q_m<0$, so $\nu>0$, see \eqref{parameternu}. Then \eqref{eq:pure-bessel-equation} becomes
\[
 \tau^2V''+\tau V'+(\tau^2-\nu^2)V=0,
\]
which is Bessel's equation of order $\nu$; see
\cite[Eq.(10.2.1)]{DLMF}.
A fundamental system is given by the Bessel functions
$J_\nu=J_\nu(z)$ and $Y_\nu=Y_\nu(z)$, $z\in\mathbb C$; see
\cite[Eqs.~(10.2.2)--(10.2.3)]{DLMF}.

\noindent For fixed $0\le s\le t$, let
\begin{align}\label{def:xylamb}
 x=e^sr,\qquad y=e^tr,\qquad
 \lambda=\frac{y}{x}=e^{t-s}.
\end{align}

\noindent The Wronskian identity
\begin{equation}\label{eq:pure-bessel-wronskian}
 J_\nu(z)Y_\nu'(z)-J_\nu'(z)Y_\nu(z)
 =\frac{2}{\pi z}
\end{equation}
follows from \cite[Eq.(10.5.2)]{DLMF}.  Solving the two corresponding
$2\times2$ systems for the initial conditions at $t=s$, we obtain
\begin{align}
 \widehat K_0(t,s,\xi)
 &=
 \frac{\pi x}{2}
 \bigl(
  Y_\nu'(x)J_\nu(y)-J_\nu'(x)Y_\nu(y)
 \bigr),
 \label{eq:pure-mult-K0}\\
 \widehat K_1(t,s,\xi)
 &=
 \frac{\pi}{2}
 \bigl(
  J_\nu(x)Y_\nu(y)-Y_\nu(x)J_\nu(y)
 \bigr).
 \label{eq:pure-mult-K1}
\end{align}
Indeed, \eqref{eq:pure-bessel-wronskian} gives
\[
 \widehat K_0(s,s,\xi)=1,\qquad
 \partial_t\widehat K_0(s,s,\xi)=0,
\]
and
\[
 \widehat K_1(s,s,\xi)=0,\qquad
 \partial_t\widehat K_1(s,s,\xi)=1.
\]
Since $\partial_t=y\partial_y$, differentiating
\eqref{eq:pure-mult-K0}--\eqref{eq:pure-mult-K1} gives
\begin{align}
 \partial_t\widehat K_0(t,s,\xi)
 &=
 \frac{\pi xy}{2}
 \bigl(
  Y_\nu'(x)J_\nu'(y)-J_\nu'(x)Y_\nu'(y)
 \bigr),
 \label{eq:pure-mult-dtK0}\\
 \partial_t\widehat K_1(t,s,\xi)
 &=
 \frac{\pi y}{2}
 \bigl(
  J_\nu(x)Y_\nu'(y)-Y_\nu(x)J_\nu'(y)
 \bigr).
 \label{eq:pure-mult-dtK1}
\end{align}
Moreover, multiplication by $y=e^tr$ gives the Fourier multiplier
corresponding to the elastic component $e^t\nabla$.  Therefore, for
$j=0,1$, it is sufficient to estimate
\[
 |\widehat K_j|,
 \qquad
 |\partial_t\widehat K_j|,
 \qquad
 y|\widehat K_j|.
\]
We divide the extended phase space into the following three regions:
\[
 \mathcal Z_1=\{0<x\le y\le1\},\qquad
 \mathcal Z_2=\{0<x\le1\le y\},\qquad
 \mathcal Z_3=\{1\le x\le y\}.
\]

\smallskip
\noindent\emph{Small arguments: $\mathcal Z_1$.}

For fixed $\nu>0$, the limiting formulas
\cite[Eqs.(10.7.3)--(10.7.4)]{DLMF},
and the recurrence identities
\[
 J_\nu'(z)=\frac{\nu}{z}J_\nu(z)-J_{\nu+1}(z),
 \qquad
 Y_\nu'(z)=\frac{\nu}{z}Y_\nu(z)-Y_{\nu+1}(z),
\]
see \cite[Eq.(10.6.2)]{DLMF},
 imply for $0<z\le1$,
\begin{equation}\label{eq:pure-bessel-small}
 \begin{aligned}
 |J_\nu(z)|&\le Cz^\nu,
 &\qquad
 |J_\nu'(z)|&\le Cz^{\nu-1},\\
 |Y_\nu(z)|&\le Cz^{-\nu},
 &
 |Y_\nu'(z)|&\le Cz^{-\nu-1}.
 \end{aligned}
\end{equation}
Here and below the constant may depend on the fixed order $\nu$, but is
independent of $x$, $y$, $s$, and $t$.

Substituting \eqref{eq:pure-bessel-small} into
\eqref{eq:pure-mult-K0}--\eqref{eq:pure-mult-dtK1} and recalling
the definitions in \eqref{def:xylamb}, we obtain
\begin{equation}\label{eq:pure-Z1}
 |\widehat K_j|
 +|\partial_t\widehat K_j|
 +y|\widehat K_j|
 \le C\left(\frac{y}{x}\right)^\nu
 =Ce^{\nu(t-s)},
 \qquad j=0,1.
\end{equation}
For clarity, the largest factor occurring in
\(
 x^{-\nu}y^\nu=\lambda^\nu.
\)
For instance,
\[
 y\,|Y_\nu(x)J_\nu'(y)|
 \lesssim yx^{-\nu}y^{\nu-1}
 =x^{-\nu}y^\nu.
\]
All other products are bounded in the same way, since
$0<x\le y\le1$.

\noindent Plancherel's theorem now yields
\begin{align}
\mathcal E_t(K_0^{\mathrm{low}}(t,s)h)
 &\le Ce^{\nu(t-s)}\norm h_2,
 \label{eq:pure-Z1-K0} \\\mathcal E_t(K_1^{\mathrm{low}}(t,s)h)
 &\le Ce^{\nu(t-s)}\norm h_2.
 \label{eq:pure-Z1-K1}
 \end{align}

\smallskip
\noindent\emph{Large arguments: $\mathcal Z_3$.}
For fixed order $\nu$, the large-argument expansions
\cite[Eqs.(10.17.3)--(10.17.4)]{DLMF} and their differentiated
versions \cite[Eqs.(10.17.9)--(10.17.10)]{DLMF} imply that
\begin{equation}\label{eq:pure-bessel-large}
 |J_\nu(z)|+|Y_\nu(z)|
 +|J_\nu'(z)|+|Y_\nu'(z)|
 \le Cz^{-\frac{1}{2}},
 \qquad z\ge1.
\end{equation}
Using \eqref{eq:pure-bessel-large} in
\eqref{eq:pure-mult-K0}--\eqref{eq:pure-mult-dtK1}, we find
\begin{align}|\widehat K_0|
 &\le C\left(\frac xy\right)^{\frac{1}{2}}
 =C\lambda^{-\frac{1}{2}},
 &
 |\partial_t\widehat K_0|+y|\widehat K_0|
 &\le C(xy)^{\frac{1}{2}}
 =Cx\lambda^{\frac{1}{2}},
 \label{eq:pure-Z3-K0-mult}
\\|\widehat K_1|
 &\le C(xy)^{-\frac{1}{2}},
 &
 |\partial_t\widehat K_1|+y|\widehat K_1|
 &\le C\left(\frac yx\right)^{\frac{1}{2}}
 =C\lambda^{\frac{1}{2}}.
 \label{eq:pure-Z3-K1-mult}
 \end{align}
Since $x=e^sr$, the factor $x$ in
\eqref{eq:pure-Z3-K0-mult} corresponds to $e^s|\nabla|$.
Consequently, Plancherel's theorem gives
\begin{align} \mathcal E_t(K_0^{\mathrm{high}}(t,s)h)
 &\le Ce^{\frac{t-s}{2}}
 \bigl(\norm h_2+e^s\norm{\nabla h}_2\bigr),
 \label{eq:pure-Z3-K0}\\
 \mathcal E_t(K_1^{\mathrm{high}}(t,s)h)
 &\le Ce^{\frac{t-s}{2}}\norm h_2.
 \label{eq:pure-Z3-K1}
 \end{align}
Thus the derivative of the initial displacement is required in the
high-frequency estimate for $K_0$.

\smallskip
\noindent\emph{Intermediate arguments: $\mathcal Z_2$.}
In this region, $0<x\le1\le y$.  We therefore apply
\eqref{eq:pure-bessel-small} to the factors evaluated at $x$ and
\eqref{eq:pure-bessel-large} to those evaluated at $y$.  From
\eqref{eq:pure-mult-K0}--\eqref{eq:pure-mult-dtK1}, we obtain
\begin{eqnarray*} |\widehat K_0|
 &&\le Cy^{-\frac{1}{2}}(x^\nu+x^{-\nu})\le Cy^{-\frac{1}{2}}x^{-\nu},
 \label{eq:pure-Z2-K0-zero}\\
 |\partial_t\widehat K_0|+y|\widehat K_0|
 &&\le Cy^{\frac{1}{2}}(x^\nu+x^{-\nu})\le Cy^{\frac{1}{2}}x^{-\nu},
 \label{eq:pure-Z2-K0-mult}\\
 |\widehat K_1|
 &&\le Cy^{-\frac{1}{2}}(x^\nu+x^{-\nu})
 \le Cy^{-\frac{1}{2}}x^{-\nu},
 \label{eq:pure-Z2-K1-zero}\\
 |\partial_t\widehat K_1|+y|\widehat K_1|
 &&\le Cy^{\frac{1}{2}}(x^\nu+x^{-\nu})
 \le Cy^{\frac{1}{2}}x^{-\nu}.
 \label{eq:pure-Z2-K1-mult}
 \end{eqnarray*}
Here we used $x^\nu\le x^{-\nu}$ for $0<x\le1$.

\noindent Since $y=e^tr$ and $x=e^sr$, the largest factor on the right-hand
side is
\begin{equation*}\label{eq:pure-transition-factor}
 y^{\frac{1}{2}}x^{-\nu}
 =e^{\frac{t}{2}-\nu s}r^{\frac{1}{2}-\nu}.
\end{equation*}
Moreover, the definition of $\mathcal Z_2$ gives
\(
 e^{-t}\le r\le e^{-s}.
\)
We distinguish two cases. \\
 If $0<\nu\le \frac{1}{2}$, then
$\frac{1}{2}-\nu\ge0$, and hence, using $r\le e^{-s}$, it follows
\[
 e^{\frac{t}{2}-\nu s}r^{\frac{1}{2}-\nu}
 \le
 e^{\frac{t}{2}-\nu s}e^{-s(\frac{1}{2}-\nu)}
 =e^{\frac{t-s}{2}}.
\]
If $\nu\ge \frac{1}{2}$, then $\frac{1}{2}-\nu\le0$, and using $r\ge e^{-t}$, it follows
\[
 e^{\frac{t}{2}-\nu s}r^{\frac{1}{2}-\nu}
 \le
 e^{\frac{t}{2}-\nu s}e^{-t(\frac{1}{2}-\nu)}
 =e^{\nu(t-s)}.
\]
Therefore,
\begin{equation*}\label{eq:pure-transition-key}
 y^{\frac{1}{2}}x^{-\nu}
 \le
 e^{\max\{\frac{1}{2},\nu\}(t-s)}.
\end{equation*}
The bounds for the moduli of $\widehat K_0$ and  $\widehat K_1$ are smaller because
\[
 y^{-\frac{1}{2}}x^{-\nu}
 =y^{-1}\bigl(y^{\frac{1}{2}}x^{-\nu}\bigr)
 \le y^{\frac{1}{2}}x^{-\nu},
\]
as $y\ge1$.  Hence, Plancherel's theorem and
\eqref{eq:pure-transition-key} yield
\begin{align} \mathcal E_t(K_0^{\mathrm{mid}}(t,s)h)
 &\le
 Ce^{\max\{\frac{1}{2},\nu\}(t-s)}
 \bigl(\norm h_2+e^s\norm{\nabla h}_2\bigr),
 \label{eq:pure-Z2-K0}\\
 \mathcal E_t(K_1^{\mathrm{mid}}(t,s)h)
 &\le
 Ce^{\max\{\frac{1}{2},\nu\}(t-s)}\norm h_2.
 \label{eq:pure-Z2-K1}
 \end{align}
Thus the derivative of the initial displacement is required in the
intermediate-frequency estimate for $K_0$. \\
Combining
\eqref{eq:pure-Z1-K1}--\eqref{eq:pure-Z1-K0},
\eqref{eq:pure-Z3-K1}--\eqref{eq:pure-Z3-K0}, and
\eqref{eq:pure-Z2-K1}--\eqref{eq:pure-Z2-K0}, we obtain
\eqref{eq:pure-K1}--\eqref{eq:pure-K0} for $q_m<0$.
\smallskip

\smallskip
\noindent\emph{The case $q_m>0$.}
Recall that $\nu=\sqrt{q_m}$ in this case. Then
\eqref{eq:pure-bessel-equation} is Bessel's equation with order
$i\nu$.  It is convenient to use the fundamental system
\[
 J_{i\nu}(z),\qquad J_{-i\nu}(z).
\]
Their Wronskian is given by \cite[Eq.(10.5.1)]{DLMF}:
\begin{equation}\label{eq:pure-imaginary-wronskian}
 \mathcal W\{J_{i\nu},J_{-i\nu}\}(z)
 =
 -\frac{2\sin(i\pi\nu)}{\pi z}
 =
 -\frac{2i\sinh(\pi\nu)}{\pi z},
\end{equation}
which does not vanish because $\nu>0$.

The power-series representation \cite[Eq.(10.2.2)]{DLMF}
and its termwise derivative give
\[
 J_{\pm i\nu}(z)
 =
\frac{1}{\Gamma(1\pm i\nu)} \Big(\frac{z}{2}\Big)^{\pm i\nu}
      \bigl(1+O(z^2)\bigr),
\]
and
\[
 J_{\pm i\nu}'(z)
 =
\frac{1}{\Gamma(1\pm i\nu)} \Big(\frac{z}{2}\Big)^{\pm i\nu}
      \left(
   \frac{\pm i\nu}{z}\bigl(1+O(z^2)\bigr)
   +O(z)
 \right).
\]
Here, for $z>0$,
\[
 \Big(\frac{z}{2}\Big)^{\pm i\nu}
 =
 e^{\pm i\nu\log\big(\frac{z}{2}\big)},
 \qquad
 \Big|\Big(\frac{z}{2}\Big)^{\pm i\nu}\Big|=1.
\]
Since $\Gamma(1\pm i\nu)\neq0$, we immediately obtain that $J_{\pm i\nu}(z)$ is bounded
near the origin and $J'_{\pm i\nu}(z)$ has potentially singular behavior of order $z^{-1}$. Therefore,
it follows that
\begin{equation}\label{eq:pure-imaginary-small}
 |J_{\pm i\nu}(z)|
 +z|J_{\pm i\nu}'(z)|
 \le C_\nu,
 \qquad 0<z\le1.
\end{equation}
Thus the  multipliers related to low arguments are uniformly bounded and carry no
factor of the form $e^{\nu(t-s)}$.

For $z\ge1$, the large argument expansions
\cite[Eqs.(10.17.3), (10.17.9)]{DLMF}, which are valid for fixed
complex order, give
\[
 |J_{\pm i\nu}(z)|
 +|J_{\pm i\nu}'(z)|
 \le C_\nu z^{-\frac{1}{2}}.
\]
Repeating the estimates in $\mathcal Z_2$ and $\mathcal Z_3$, with
\eqref{eq:pure-imaginary-small} in place of
\eqref{eq:pure-bessel-small}, yields
\[
 \mathcal E_t(K_0(t,s)h)
 \le Ce^{\frac{t-s}{2}}
 \bigl(\norm h_2+e^s\norm{\nabla h}_2\bigr),
\]
and
\[
 \mathcal E_t(K_1(t,s)h)
 \le Ce^{\frac{t-s}{2}}\norm h_2.
\]

\smallskip
\noindent\emph{The case $q_m=0$.}
In this case a fundamental system is $J_0,Y_0$.  The limiting formulas
\cite[Eqs.(10.7.1)--(10.7.2)]{DLMF} give
\[
 J_0(z)=1+O(z^2),
 \qquad
 Y_0(z)=\frac{2}{\pi}\log z+O(1)
 \qquad (z\downarrow0).
\]
Moreover, by \cite[Eq.(10.6.3)]{DLMF} and
\cite[Eqs.(10.7.3)--(10.7.4)]{DLMF},
\[
 |J_0'(z)|\le Cz,
 \qquad
 |Y_0'(z)|\le Cz^{-1},
 \qquad 0<z\le1.
\]
When $0<x\le y\le1$, the logarithmic terms occur through
differences of the form
\(
 \log y-\log x=t-s.
\)
Consequently, the low arguments estimates acquire at most the factor
\[
 1+t-s=\Omega(t-s).
\]
The large arguments estimates are unchanged, so the high and
intermediate frequency contributions are bounded by $Ce^{\frac{t-s}{2}}$.
Since, for every $\delta>0$,
\(
 1+t-s\le C_\delta e^{\delta(t-s)},
\)
this proves the homogeneous estimates for every value of $q_m$, with
the borderline factor $\Omega$ stated in the lemma.

Finally, if $q_m\neq0$, Minkowski's integral inequality and
\eqref{eq:pure-K1} give
\begin{align*}
 e^{-\beta_0 t}\mathcal E_t\left(\int_0^tK_1(t,s)F(s)\,ds\right)
 &\le
 C\int_0^t
 e^{-\beta_0 t}e^{\beta_0(t-s)}
 \norm{F(s)}_2\,ds\\
 &\le C\int_0^T e^{-\beta_0 s}\norm{F(s)}_2\,ds.
\end{align*}
Taking the supremum over $0\le t\le T$ proves
\eqref{eq:pure-duhamel}. If $q_m=0$, we use
$1+t-s\le C_\delta e^{\delta(t-s)}$ in the same calculation, with
$\beta_0+\delta$ in place of $\beta_0$, and obtain
\eqref{eq:pure-duhamel-borderline}. This completes the proof.
\end{proof}

After these estimates for solutions of the family of linear parameter-dependent Cauchy problems are established, we are now in a position to prove
 \Cref{thm:global}. The argument combines the propagator
bounds with the corresponding Duhamel estimates and a standard fixed-point
procedure in a suitable exponentially weighted energy space.

\begin{proof}[Proof of \Cref{thm:global}]
For $T>0$, we introduce the exponentially weighted energy space
\[
 X(T):=C\big([0,T],H^1(\R^d)\big)
 \cap C^1\big([0,T],L^2(\R^d)\big),
\]
endowed with the norm
\[
 \norm{u}_{X(T)}
 :=
 \sup_{0\le t\le T}e^{-\beta t}
 \left(
  \norm{u(t)}_{L^2}
  +\norm{u_t(t)}_{L^2}
  +e^t\norm{\nabla u(t)}_{L^2}
 \right).
\]
Fix
\[
 \beta\in\Big(\beta_0,\frac{n+d}{2}\Big).
\]
This interval is nonempty. Indeed, if \(q_m<0\), then
\(\beta_0=\nu\), and
\[
 m^2>-\frac{d(2n+d)}{4}
 \quad\Longleftrightarrow\quad
 \nu<\frac{n+d}{2};
\]
see \eqref{RestmTheo1}. If \(q_m\geq0\), then
\(\beta_0=\frac12<\frac{n+d}{2}\); see \eqref{def:beta}.
In the borderline case \(q_m=0\), we apply
\eqref{eq:pure-duhamel-borderline} with
\[
 \delta=\beta-\beta_0>0.
\]

\noindent The Duhamel map $\mathcal N$ for \eqref{eq:u} is
\[
\mathcal N: u(t) \to  \mathcal Nu(t)=K_0(t,0)u_0+K_1(t,0)u_1
 +\int_0^tK_1(t,s)e^{-\frac n2(p-1)s}|u(s)|^p\,ds.
\]
By \Cref{lem:pure-L2}, its homogeneous part is bounded by $C_\beta e^{\beta T}(\norm{u_0}_{H^1}+\norm{u_1}_2)$.

\noindent Put $\theta=d\frac{p-1}{2p}$.  Condition \eqref{eq:p-range} is exactly
$0<\theta\le1$.  Gagliardo--Nirenberg inequality gives
\begin{align}
 \norm{|u(s)|^p}_2
 &=\norm{u(s)}_{2p}^{p}\notag\le C\norm{u(s)}_2^{p(1-\theta)}
          \norm{\nabla u(s)}_2^{p\theta}\notag\\
 &\le C e^{(p\beta-\frac d2(p-1))s}\norm u_{X(T)}^p.
 \label{eq:pure-source}
\end{align}
Multiplying by the nonlinear coefficient and by the Duhamel weight gives
\begin{align*}
 e^{-\beta s}e^{-\frac n2(p-1)s}\norm{|u(s)|^p}_2
 &\le C e^{-(p-1)\frac{n+d}{2}-\beta)s}\norm u_{X(T)}^p.
\end{align*}
The exponential term is integrable over $[0,\infty)$ because
$\beta<\frac{n+d}{2}$.  Hence, \eqref{eq:pure-duhamel} implies, uniformly in $T$,
\[
 \norm{\mathcal Nu}_{X(T)}
 \le C\varepsilon+C\norm u_{X(T)}^p.
\]
For the difference, the pointwise inequality
\[
 \bigl||u|^p-|v|^p\bigr|
 \le C_p(|u|^{p-1}+|v|^{p-1})|u-v|
\]
and H\"older's inequality imply
\[
 \norm{|u|^p-|v|^p}_2
 \le C\bigl(\norm u_{2p}^{p-1}+\norm v_{2p}^{p-1}\bigr)
       \norm{u-v}_{2p}.
\]
Applying the same Gagliardo--Nirenberg estimate to every factor yields
\[
 \norm{\mathcal Nu-\mathcal Nv}_{X(T)}
 \le C\bigl(\norm u_{X(T)}^{p-1}+\norm v_{X(T)}^{p-1}\bigr)
       \norm{u-v}_{X(T)}.
\]
Choose $R=2C\varepsilon$.  If $\varepsilon$ is sufficiently small,
$\mathcal N$ maps the closed ball of radius $R$ in $X(T)$ into itself and
is a contraction there, with constants independent of $T$.  Banach's
fixed-point theorem gives a solution on every $[0,T]$, and uniqueness
makes these solutions compatible as $T$ varies.  Thus $u$, and hence,
$\phi=e^{-n \frac{t}{2}}u$, is the unique global (in time) energy solution.  Moreover, the
fixed-point estimate gives
\[
 \norm{u(t)}_2+\norm{u_t(t)}_2
 +e^t\norm{\nabla u(t)}_2
 \le C_\beta\varepsilon e^{\beta t}.
\]
Since due to \eqref{eq:transform} it holds
\[
 \phi=e^{-n\frac{t}{2}}u,\qquad
 \phi_t=e^{-n\frac{t}{2}}\left(u_t-\frac n2u\right),\qquad
 \nabla\phi=e^{-n\frac{t}{2}}\nabla u,
\]
the preceding inequality immediately yields
\eqref{eq:phi-energy-estimate}.
\end{proof}

%%%%%%%%
%%%%%%%%
%%%%%%%%

\subsection{Proof of \Cref{thm:global-weighted}} \label{Section3.2}
The proof is based on a reduction that is particularly suited to the
intermediate mass regime. After passing to the conformal time, we re-scale
the unknown in such a way that the mass term disappears and the transformed
equation acquires a scale-invariant damping term. A direct argument in the
non-weighted energy space, along the lines of the proof of
\Cref{thm:global}, would require \(p\geq2\) and would therefore fail to
cover the full range \(p>1\) considered here. To overcome this restriction,
we exploit the dissipative structure of the transformed equation in an
exponentially weighted energy space whose Gaussian weight is adapted to
its propagation. This framework, inspired by
\cite{IkehataTanizawa,Palmieri2018}, provides the weighted control of the
nonlinearity needed when \(1<p<2\). It allows us to combine estimates for solutions to
the associated family of parameter-dependent linear Cauchy problems with
the nonlinear energy estimates required in the fixed-point argument.

Set
\begin{equation}\label{eq:change}
 \tau=e^t-1 \quad \mbox{and} \quad
 \phi(t,x)=(1+\tau)^\rho v(\tau,x),
\end{equation}
where $\rho$ is defined in \eqref{eq:rho}
Then the Cauchy problem \eqref{eq:phi} is equivalent to
\begin{equation}\label{eq:damped-v}
 \begin{cases}
 v_{\tau\tau}-\Delta v+\dfrac{\mu}{1+\tau}v_\tau
 =(1+\tau)^\alpha|v|^p,\\[2mm]
 v(0,x)=f(x),\qquad v_\tau(0,x)=g(x)-\rho f(x),
 \end{cases}
\end{equation}
where
\begin{equation}\label{eq:mu-alpha}
 \mu=n+1+2\rho,\qquad \alpha=\rho(p-1)-2.
\end{equation}
Notice that, by \eqref{eq:mu-alpha} and the definition of $\rho$,
\[
 \frac{\mu-1}{2}
 =
 \frac n2+\rho
 =
 \sqrt{\frac{n^2}{4}-m^2}
 =
 \nu,
\]
where $\nu$ is the parameter introduced in \eqref{parameternu}. Hence, the
order of the Bessel equation arising below coincides with the parameter
already used in the proof of \Cref{thm:global}.
\begin{remark}\label{rem:rho}
The mass condition \eqref{rest:thm2} in \Cref{thm:global-weighted} is equivalently
expressed as
\[
 \frac d2\leq\rho<d,
\]
where the first inequality is strict when \(n=1\). Consequently,
\[
 \mu=n+1+2\rho>d+2.
\]
Indeed, if \(n>1\), then
\[
 \mu\geq n+d+1>d+2,
\]
whereas, if \(n=1\), the strict inequality \(\rho>\frac{d}{2}\) gives
\(\mu=2+2\rho>d+2\). Thus, \eqref{eq:damped-v} lies in the strictly
effective damped regime required by \Cref{lem:linear} below.
\end{remark}
Following \cite{IkehataTanizawa,Palmieri2018}, we introduce the Gaussian weight
\begin{equation}\label{eq:Psi}
 \Psi(\tau,x):=\frac{\mu|x|^2}{2(1+\tau)^2},
\end{equation}
and, for a pair $(h_0,h_1)$, set
\begin{equation*}
 \norm{(h_0,h_1)}_{\mathcal A_\mu}
 :=\norm{e^{\Psi(0)}h_0}_{H^1}
 +\norm{e^{\Psi(0)}h_1}_2,\,\,\,\Psi(0)=\Psi(0,x).
\end{equation*}

\begin{remark}\label{l1estimatives}
The Gaussian weight \eqref{eq:Psi}
is adapted to both the spatial scale \(1+\tau\) and the damping
coefficient in \eqref{eq:damped-v}. Indeed, it satisfies
\begin{align}\label{Rem33deriva}
 \Psi_\tau=-\frac{2}{1+\tau}\Psi,
 \qquad
 |\nabla\Psi|^2+\frac{\mu}{1+\tau}\Psi_\tau=0.
\end{align}
The second identity provides the exact cancelation required in the
weighted energy estimate.
Moreover, the Gaussian decay associated with \(\Psi\) yields the
\(L^1\)-integrability needed in the desired
estimates for solutions of a family of parameter-dependent Cauchy problems. More precisely,
for every \(\tau\geq0\),
\begin{align*}
 \norm{e^{-\Psi(\tau)}}_{L^2}^2
 =
 \int_{\R^d}
 \exp\left(-\frac{\mu|x|^2}{(1+\tau)^2}\right)\dd x =
 \left(\frac{\pi}{\mu}\right)^{\frac{d}{2}}(1+\tau)^d.
\end{align*}
In particular, setting \(\tau=0\), we obtain
\[
 \norm{e^{-\Psi(0)}}_{L^2}
 =
 \left(\frac{\pi}{\mu}\right)^{\frac{d}{4}}.
\]
Hence, by the Cauchy--Schwarz inequality,
\[
 \norm{h}_{L^1}
 =
 \norm{e^{-\Psi(0)}e^{\Psi(0)}h}_{L^1}
 \leq
 \left(\frac{\pi}{\mu}\right)^{\frac{d}{4}}
 \norm{e^{\Psi(0)}h}_{L^2}.
\]
Consequently, the weighted assumption
\(e^{\Psi(0)}h\in L^2(\R^d)\) implies that \(h\in L^1(\R^d)\).
\end{remark}

For $\tau\ge s\ge0$, let $\K_0=\K_0(\tau,s)$ and $\K_1=\K_1(\tau,s)$ be the
propagators associated with the linear Cauchy problem \eqref{eq:lineardamped-v} below.
We also define
\begin{equation}\label{eq:damped-energy}
 \E_\tau(w):=\norm{w(\tau)}_2
 +(1+\tau)\norm{(\nabla w,w_\tau)(\tau)}_2.
\end{equation}
The following estimates are obtained by specializing
\cite[Propositions~5.1 and~5.2]{Palmieri2018} to the massless case.
For the sake of completeness, we provide the proof in the present setting.
\begin{lemma}\label{lem:linear}
Assume that $\mu> d+2$. If $v$ solves the Cauchy problem \begin{equation}\label{eq:lineardamped-v}
 v_{\tau\tau}-\Delta v+\dfrac{\mu}{1+\tau}v_\tau
 =0, \quad
 v(0,x)=f(x),\quad v_\tau(0,x)=g(x)-\rho f(x),
 \end{equation} then
\begin{equation}\label{eq:hom}
 \mathcal E_\tau(v)
 \leq
 C(1+\tau)^{-\frac{d}{2}}
 \norm{(v(0),v_\tau(0))}_{\mathcal A_\mu}.
\end{equation}
For the solution of the Cauchy problem \begin{equation}\label{eq:sourcelineardamped-v}
 v_{\tau\tau}-\Delta v+\dfrac{\mu}{1+\tau}v_\tau
 =F(\tau), \quad
 v(0)=0,\quad v_\tau(0)=0,
 \end{equation} it follows that
\begin{equation}\label{eq:duhamel}
 \mathcal E_\tau\Big(\int_0^\tau\K_1(\tau,s)F(s)\,\dd s\Big)
 \leq
 C(1+\tau)^{-\frac{d}{2}}
 \int_0^\tau
 \Big(
  (1+s)\norm{F(s)}_{L^1}
  +(1+s)^{1+\frac{d}{2}}\norm{F(s)}_{L^2}
 \Big)\dd s.
\end{equation}
\end{lemma}

\begin{proof}
Taking the Fourier transform with respect to $x$ in \eqref{eq:lineardamped-v} becomes
\[
 \widehat w_{\tau\tau}
 +\frac{\mu}{1+\tau}\widehat w_\tau
 +|\xi|^2\widehat w
 =0.
\]
Set
\[
 \nu:=\frac{\mu-1}{2}>\frac{d+1}{2},
 \qquad
 r:=(1+\tau)|\xi|,
 \qquad
 \widehat w(\tau,\xi):=(1+\tau)^{-\nu}Z(r).
\]
A direct computation shows that $Z$ satisfies Bessel's equation
\begin{equation}\label{eq:bessel}
 r^2Z_{rr}+rZ_r+(r^2-\nu^2)Z=0.
\end{equation}
For fixed $0\leq s\leq\tau$, we introduce the notation
\[
 a:=1+s,\qquad
 b:=1+\tau,\qquad
 x:=a|\xi|,\qquad
 y:=b|\xi|.
\]
In particular,
\[
 1\leq a\leq b,
 \qquad
 0<x\leq y,
 \qquad
 \frac{x}{y}=\frac{a}{b}.
\]
Proceeding as in the proof of  \Cref{lem:pure-L2}, and solving the
corresponding systems for the initial conditions at $\tau=s$, we obtain
\begin{align}
 \widehat\K_0(\tau,s,\xi)
 &=
 \frac{\pi x}{2}\left(\frac ab\right)^\nu
 \bigl(
  J_{\nu+1}(x)Y_\nu(y)
  -Y_{\nu+1}(x)J_\nu(y)
 \bigr),
 \label{eq:K0}\\
 \widehat\K_1(\tau,s,\xi)
 &=
 \frac{\pi a}{2}\left(\frac ab\right)^\nu
 \bigl(
  J_\nu(x)Y_\nu(y)
  -Y_\nu(x)J_\nu(y)
 \bigr).
 \label{eq:K1}
\end{align}
These representations follow from standard Bessel identities
\cite[Eqs.(10.5.2), (10.6.6)]{DLMF}.

To estimate the time-derivative component of the energy, we differentiate
\eqref{eq:K0} and \eqref{eq:K1}. Using
\[
 \frac{\dd}{\dd z}\bigl(z^{-\nu}J_\nu(z)\bigr)
 =
 -z^{-\nu}J_{\nu+1}(z),
 \qquad
 \frac{\dd}{\dd z}\bigl(z^{-\nu}Y_\nu(z)\bigr)
 =
 -z^{-\nu}Y_{\nu+1}(z),
\]
we obtain
\begin{align}
 \partial_\tau\widehat\K_0(\tau,s,\xi)
 &=
 -\frac{\pi x|\xi|}{2}\left(\frac ab\right)^\nu
 \bigl(
  J_{\nu+1}(x)Y_{\nu+1}(y)
  -Y_{\nu+1}(x)J_{\nu+1}(y)
 \bigr),
 \label{eq:dtK0}\\
 \partial_\tau\widehat\K_1(\tau,s,\xi)
 &=
 -\frac{\pi a|\xi|}{2}\left(\frac ab\right)^\nu
 \bigl(
  J_\nu(x)Y_{\nu+1}(y)
  -Y_\nu(x)J_{\nu+1}(y)
 \bigr).
 \label{eq:dtK1}
\end{align}
Thus, the solution, spatial-derivative, and time-derivative components
of the energy can all be estimated from
\eqref{eq:K0}--\eqref{eq:dtK1}.

We divide the extended phase space according to the sizes of $x$ and $y$ as follows:
\[
 Z_1:=\{0<x\leq y\leq1\},\qquad
 Z_2:=\{0<x\leq1\leq y\},\qquad
 Z_3:=\{1\leq x\leq y\}.
\]
For small arguments, standard asymptotic formulas and recurrence
identities yield
\[
 |J_\gamma(z)|\lesssim z^\gamma,
 \qquad
 |J_\gamma'(z)|\lesssim z^{\gamma-1},
 \qquad
 |Y_\gamma(z)|\lesssim z^{-\gamma},
 \qquad
 |Y_\gamma'(z)|\lesssim z^{-\gamma-1},
 \qquad 0<z\leq1,
\]
for $\gamma=\nu,\nu+1$. For large arguments, we use
\[
 |J_\gamma(z)|+|Y_\gamma(z)|
 +|J_\gamma'(z)|+|Y_\gamma'(z)|
 \lesssim z^{-\frac{1}{2}},
 \qquad z\geq1.
\]
We now illustrate how these pointwise estimates lead to the desired
operator bounds. We consider $\K_1(\tau,s)h$ only; the other components are
treated similarly.

 First, we turn to $Z_1$. Since $0<x\leq y\leq1$, the small argument
bounds applied to \eqref{eq:K1} give
\begin{align*}
 |\widehat\K_1(\tau,s,\xi)|
 \lesssim
 a\Big(\frac ab\Big)^\nu
 \bigl(x^\nu y^{-\nu}+x^{-\nu}y^\nu\bigr)=
 a\Big(
  \Big(\frac{x}{y}\Big)^{2\nu}+1
 \Big)
 \lesssim a.
\end{align*}
Moreover, since
 \(
 Z_1\subset\{\,|\xi|\leq b^{-1}\,\}
 \), we have \(
 |Z_1|\lesssim b^{-d}.\)
It follows that
\begin{align*}
 \big\|\widehat\K_1(\tau,s,\xi)\widehat h(\xi)\big\|_{L^2(Z_1)}
 \leq
 \sup_{\xi\in Z_1}\big|\widehat\K_1(\tau,s,\xi)\big|
 |Z_1|^{\frac{1}{2}}\|\widehat h\|_{L^\infty} \lesssim
 ab^{-\frac{d}{2}}\|h\|_{L^1}.
\end{align*}

We next consider $Z_2\cup Z_3$. Substituting the appropriate small and
large argument bounds into \eqref{eq:K1} gives
\[
 \big|\widehat\K_1(\tau,s,\xi)\big|
 \lesssim
 a\left(\frac ab\right)^{\frac{d}{2}},
 \,\,\,\mbox{in}\,\,\,Z_2\cup Z_3.
\]
More precisely, after multiplication by the factor $1+b|\xi|$ associated
with the energy, the Bessel estimates produce the power
$(\frac{a}{b})^{\nu-\frac{1}{2}}$. Since
\(
 \nu-\frac{1}{2}> \frac{d}{2}
 \,\,\,\text{and}\,\,\,
 0<\frac{a}{b}\leq1,
\)
this power is bounded by the displayed rate. Therefore, by
Plancherel's theorem,
\begin{align*}
 \big\|\widehat\K_1(\tau,s,\xi)\widehat h(\xi)\big\|_{L^2(Z_2\cup Z_3)}
 \lesssim
 a\Big(\frac ab\Big)^{\frac{d}{2}}
 \|\widehat h\|_{L^2(\mathbb R^d)} \lesssim
 a\Big(\frac ab\Big)^{\frac{d}{2}}
 \|h\|_{L^2(\mathbb R^d)}.
\end{align*}
Combining the estimates in $Z_1$ and $Z_2\cup Z_3$, we arrive at
\[
 \|\K_1(\tau,s)h\|_{L^2}
 \lesssim
 ab^{-\frac{d}{2}}\|h\|_{L^1}
 +
 a\Big(\frac ab\Big)^{\frac{d}{2}}\|h\|_{L^2}.
\]
The estimates for $\partial_\tau\K_1(\tau,s)h$ and
$\nabla\K_1(\tau,s)h$ follow from the same decomposition, now using
\eqref{eq:dtK1} and $|\xi|\widehat\K_1$, respectively. Applying the
same argument to \eqref{eq:K0} and \eqref{eq:dtK0} gives the
corresponding bounds for $\K_0(\tau,s)h$. Combining all these estimates
according to the definition of $\E_\tau$, we obtain
\begin{align}
 \E_\tau(\K_0(\tau,s)h)
 &\leq
 Cb^{-\frac{d}{2}}\|h\|_{L^1}
 +
 C\Big(\frac ab\Big)^{\frac{d}{2}}
 \big(
  \|h\|_{L^2}
  +a\|\nabla h\|_{L^2}
 \big),
 \label{eq:K0-final}\\
 \E_\tau(\K_1(\tau,s)h)
 &\leq
 Cab^{-\frac{d}{2}}\|h\|_{L^1}
 +
 Ca\Big(\frac ab\Big)^{\frac{d}{2}}\|h\|_{L^2}.
 \label{eq:K1-final}
\end{align}
The additional term $a\|\nabla h\|_{L^2}$ in
\eqref{eq:K0-final} arises from the spatial-derivative component of the
energy. Indeed, the corresponding Fourier estimate contains the factor
$1+a|\xi|$.

We now apply these estimates to \eqref{eq:lineardamped-v}.
At $s=0$, we
have
\(
 a=1.
\)
Since
\[
 w(\tau)
 =
 \K_0(\tau,0)w_0+\K_1(\tau,0)w_1,
\]
the triangle inequality, \eqref{eq:K0-final}, and
\eqref{eq:K1-final} imply
\begin{align*}
 \E_\tau(w)
 &\leq
 C(1+\tau)^{-\frac{d}{2}}
 \bigl(
  \|w_0\|_{L^1}
  +\|w_0\|_{L^2}
  +\|\nabla w_0\|_{L^2}
  +\|w_1\|_{L^1}
  +\|w_1\|_{L^2}
 \bigr).
\end{align*}
By  \Cref{l1estimatives}, the $L^1$-norms are controlled by the
corresponding Gaussian-weighted $L^2$-norms. Consequently, we proved \eqref{eq:hom}.

Finally, consider the Cauchy problem \eqref{eq:sourcelineardamped-v}.
Minkowski's integral inequality and \eqref{eq:K1-final} imply
\begin{align*}
 \int_0^\tau
 \E_\tau\bigl(\K_1(\tau,s)F(s)\bigr)\,\dd s
 &\leq
 C(1+\tau)^{-\frac{d}{2}}
 \int_0^\tau
 \left(
  (1+s)\|F(s)\|_{L^1}
  +(1+s)^{1+\frac{d}{2}}\|F(s)\|_{L^2}
 \right)\dd s.
\end{align*}
This proves \eqref{eq:duhamel} and completes the proof.
\end{proof}
\begin{remark}
At the endpoint \(\mu=d+2\), the estimates in
\Cref{lem:linear} remain valid with the additional logarithmic loss
\[
 \ell(\tau):=1+\bigl(\log(1+\tau)\bigr)^{\frac{1}{2}}
\]
on their right-hand sides. More precisely, the factor
\((1+\tau)^{-\frac{d}{2}}\) in \eqref{eq:hom} and \eqref{eq:duhamel} must be
replaced by
\[
 (1+\tau)^{-\frac{d}{2}}\ell(\tau).
\]
\end{remark}
We shall use the following weighted Gagliardo--Nirenberg inequality,
which is a direct adaptation of \cite[Lemma~4.2]{Palmieri2018} to the
weight \(\Psi\) considered here, which differs from the weight used therein
only by a multiplicative constant.

For each \(\tau\geq0\), we set
\[
 H^1_{\Psi(\tau)}
 :=
 \left\{
 h\in H^1(\R^d):
 e^{\Psi(\tau)}h,\,
 e^{\Psi(\tau)}\nabla h\in L^2(\R^d)
 \right\}.
\]
\begin{lemma}
\label{lem:GN}
Let \(2\leq q\leq \frac{2d}{(d-2)_+}\), \(0<\delta\leq1\), and
\[
 \vartheta(q):=d\left(\frac12-\frac1q\right).
\]
Then, for every \(h\in H^1_{\Psi(\tau)}(\R^d)\),
\begin{equation}\label{eq:GN}
 \norm{e^{\delta\Psi(\tau)}h}_q
 \leq
 C(1+\tau)^{1-\vartheta(q)}
 \norm{\nabla h}_2^{1-\delta}
 \norm{e^{\Psi(\tau)}\nabla h}_2^\delta.
\end{equation}
\end{lemma}

For $T>0$, introduce
\begin{align} \label{defEspX22}
 X(T):=C\big([0,T],H^1(\R^d)\big)\cap C^1\big([0,T],L^2(\R^d)\big),
\end{align}
with the norm
\begin{align}\label{eq:X-norm}
 \norm v_{X(T)}
 &:=
 \sup_{0\le\tau\le T}\left[
 (1+\tau)^{\frac{d}{2}}\norm{v(\tau)}_2
 +(1+\tau)^{\frac{d}{2}+1}\norm{(\nabla v,v_\tau)(\tau)}_2\right]\notag\\
 &\quad+
 \sup_{0\le\tau\le T}\left[
 (1+\tau)^{-1}\norm{e^{\Psi(\tau)}v(\tau)}_2
 +\norm{e^{\Psi(\tau)}(\nabla v,v_\tau)(\tau)}_2\right].
\end{align}
The next lemma provides desired estimates required in the
fixed-point argument. Its proof adapts the weighted interpolation
procedure employed in \cite{Palmieri2018} to the time-dependent
nonlinearity and the decay structure of the present problem.
\begin{lemma}
\label{lem:nonlinear}
Let $p>1$ satisfy \eqref{eq:p-range}, let $\alpha$ be defined by
\eqref{eq:mu-alpha}, and let $v\in X(T)$. Then, for every sufficiently
small $\eta>0$,
\begin{align}
\label{eq:nonlinear-source}
 &(1+\tau)\left(
 \norm{(1+\tau)^\alpha|v(\tau)|^p}_1
 +(1+\tau)^{\frac{d}{2}}
 \norm{(1+\tau)^\alpha|v(\tau)|^p}_2
 \right)\notag\\
 &\qquad\leq
 C(1+\tau)^{-1-(d-\rho)(p-1)+\eta}
 \norm v_{X(T)}^p,
 \qquad 0\leq\tau\leq T.
\end{align}
Moreover, for every $v,w\in X(T)$,
\begin{align}
\label{eq:nonlinear-source-difference}
 &(1+\tau)\left(
 \norm{(1+\tau)^\alpha\bigl(|v(\tau)|^p-|w(\tau)|^p\bigr)}_1
 \right.\left.+(1+\tau)^{\frac{d}{2}}
 \norm{(1+\tau)^\alpha\bigl(|v(\tau)|^p-|w(\tau)|^p\bigr)}_2
 \right)\notag\\
 &\qquad\leq
 C(1+\tau)^{-1-(d-\rho)(p-1)+\eta}
 \bigl(\norm v_{X(T)}^{p-1}+\norm w_{X(T)}^{p-1}\bigr)
 \norm{v-w}_{X(T)}.
\end{align}
\end{lemma}

\begin{proof}
We first prove \eqref{eq:nonlinear-source}. For every $\delta>0$,
 \Cref{l1estimatives} gives
\[
 \norm h_1\leq
 C_\delta(1+\tau)^{\frac{d}{2}}\norm{e^{\delta\Psi(\tau)}h}_2.
\]
Moreover, since $e^{\delta\Psi}\geq1$, it follows
\[
 (1+\tau)^{\frac{d}{2}}\norm h_2
 \leq(1+\tau)^{\frac{d}{2}}\norm{e^{\delta\Psi(\tau)}h}_2.
\]
Consequently,
\begin{equation}
\label{eq:gaussian-source}
 \norm h_1+(1+\tau)^{\frac{d}{2}}\norm h_2
 \leq C_\delta(1+\tau)^{\frac{d}{2}}
 \norm{e^{\delta\Psi(\tau)}h}_2.
\end{equation}
We apply \eqref{eq:gaussian-source} to
$h=(1+\tau)^\alpha|v(\tau)|^p$. After replacing $\delta$ by
$p\delta$, we obtain
\begin{align}
\label{eq:source-after-gaussian}
 &(1+\tau)\left(
 \norm{(1+\tau)^\alpha|v|^p}_1
 +(1+\tau)^{\frac{d}{2}}\norm{(1+\tau)^\alpha|v|^p}_2
 \right)\notag\\
 &\qquad\leq
 C_\delta(1+\tau)^{1+\alpha+\frac{d}{2}}
 \norm{e^{\delta\Psi(\tau)}v}_{2p}^p.
\end{align}
The restriction \eqref{eq:p-range} guarantees that $2p$ lies in the
admissible range of \Cref{lem:GN}. Applying that lemma with $q=2p$ gives
\[
 \norm{e^{\delta\Psi(\tau)}v}_{2p}^p
 \leq
 C(1+\tau)^{p(1-\vartheta(2p))}
 \norm{\nabla v}_2^{p(1-\delta)}
 \norm{e^{\Psi(\tau)}\nabla v}_2^{p\delta}.
\]
By the definition of $X(T)$ we have
\[
 \norm{\nabla v(\tau)}_2
 \leq(1+\tau)^{-\frac{d}{2}-1}\norm v_{X(T)},
 \quad
 \norm{e^{\Psi(\tau)}\nabla v(\tau)}_2
 \leq\norm v_{X(T)}.
\]
Substitution in \eqref{eq:source-after-gaussian} yields
\[
 (1+\tau)\left(
 \norm{(1+\tau)^\alpha|v|^p}_1
 +(1+\tau)^{\frac{d}{2}}\norm{(1+\tau)^\alpha|v|^p}_2
 \right)
 \leq C(1+\tau)^\Gamma\norm v_{X(T)}^p,
\]
where
\[
 \Gamma
 =1+\alpha+\frac d2+p\bigl(1-\vartheta(2p)\bigr)
 -p(1-\delta)\Big(\frac d2+1\Big).
\]
Since
\[
 \vartheta(2p)=\frac{d(p-1)}{2p},
 \qquad
 \alpha=\rho(p-1)-2,
\]
we have
\[
 \Gamma=-1-(d-\rho)(p-1)
 +p\delta\Big(\frac d2+1\Big).
\]
Choosing $\delta>0$ sufficiently small and writing
\(
 \eta:=p\delta\left(\frac d2+1\right)
\)
proves \eqref{eq:nonlinear-source}.

Finally,
\[
 \bigl||a|^p-|b|^p\bigr|
 \leq C_p\bigl(|a|^{p-1}+|b|^{p-1}\bigr)|a-b|.
\]
Applying \eqref{eq:gaussian-source}, H\"older's inequality, and
 \Cref{lem:GN} exactly as above proves
\eqref{eq:nonlinear-source-difference}.
\end{proof}

\medskip

\begin{proof}[Proof of \Cref{thm:global-weighted}]
Let us consider $X(T)$ as in \eqref{defEspX22}. For $v\in X(T)$, define the mapping
\begin{equation}
\label{eq:N-map}
\N: v \to \N v(\tau)
 :=
 \K_0(\tau,0)f+\K_1(\tau,0)(g-\rho f)
 +\int_0^\tau
 \K_1(\tau,s)(1+s)^\alpha|v(s)|^p\dd s.
\end{equation}
Thus, $u:=\N v$ is the unique energy solution of
\begin{equation}
\label{eq:linear-problem-Nv}
 \begin{cases}
 u_{\tau\tau}-\Delta u+\dfrac{\mu}{1+\tau}u_\tau
 =(1+\tau)^\alpha|v|^p,\\[2mm]
 u(0,x)=f(x),\qquad u_\tau(0,x)=g(x)-\rho f(x).
 \end{cases}
\end{equation}
% Notice that $v$ is an arbitrary element of $X(T)$; no solution of the
% nonlinear problem is being presupposed. Only after a fixed point of
% $\N$ has been obtained will \eqref{eq:linear-problem-Nv} become the
% nonlinear equation.

We estimate separately the entries of
\eqref{eq:X-norm}. From  \Cref{lem:linear} and
\eqref{eq:nonlinear-source}, we obtain
\begin{align}
\label{eq:N-unweighted}
 &\sup_{0\leq\tau\leq T}\left[
 (1+\tau)^{\frac{d}{2}}\norm{\N v(\tau)}_2
 +(1+\tau)^{\frac{d}{2}+1}
 \norm{(\nabla\N v,\partial_\tau\N v)(\tau)}_2
 \right]\notag\\
 &\qquad\leq
 C\eps+C\norm v_{X(T)}^p
 \int_0^T(1+s)^{-1-(d-\rho)(p-1)+\eta}\dd s.
\end{align}
We choose
\begin{equation}
\label{eq:eta-choice}
 0<\eta<(d-\rho)(p-1).
\end{equation}
Since $\rho<d$ and \eqref{eq:eta-choice}, the time integral in \eqref{eq:N-unweighted} is
uniformly bounded, independently of $T$. Hence,
\begin{align}
\label{eq:N-unweighted-final}
 &\sup_{0\leq\tau\leq T}\left[
 (1+\tau)^{\frac{d}{2}}\norm{\N v(\tau)}_2
 +(1+\tau)^{\frac{d}{2}+1}
 \norm{(\nabla\N v,\partial_\tau\N v)(\tau)}_2
 \right]\notag\\
 &\qquad\leq C\eps+C\norm v_{X(T)}^p.
\end{align}
We next estimate the exponentially weighted component. Set
\[
 \E_\Psi[u](\tau)
 :=\frac12\int_{\R^d}e^{2\Psi(\tau,x)}
 \bigl(|u_\tau(\tau,x)|^2+|\nabla u(\tau,x)|^2\bigr)\dd x.
\]
Multiplying \eqref{eq:linear-problem-Nv} by $e^{2\Psi}u_\tau$ and
integrating over $\R^d$, we obtain
\begin{align}
\label{eq:weighted-energy-identity}
 \frac{\dd}{\dd\tau}\E_\Psi[u](\tau)
 &=\int_{\R^d}e^{2\Psi}(1+\tau)^\alpha|v|^p u_\tau\dd x
 -\frac{\mu}{1+\tau}\int_{\R^d}e^{2\Psi}|u_\tau|^2\dd x\notag\\
 &\quad+\int_{\R^d}\Psi_\tau e^{2\Psi}
 \bigl(|u_\tau|^2+|\nabla u|^2\bigr)\dd x
 -2\int_{\R^d}e^{2\Psi}u_\tau\nabla\Psi\cdot\nabla u\dd x.
\end{align}
Using \eqref{Rem33deriva}
and Young's inequality, we have
\begin{align*}
 -2u_\tau\nabla\Psi\cdot\nabla u
 &\leq
 \frac{\mu}{1+\tau}|u_\tau|^2
 +\frac{1+\tau}{\mu}|\nabla\Psi|^2|\nabla u|^2=\frac{\mu}{1+\tau}|u_\tau|^2
 -\Psi_\tau|\nabla u|^2.
\end{align*}
Therefore, the sum of the last three terms on the right-hand side of
\eqref{eq:weighted-energy-identity} is non-positive. Hence,
\begin{align}
\label{eq:weighted-energy-differential}
 \frac{\dd}{\dd\tau}\E_\Psi[u](\tau)
 &\leq
 \norm{e^{\Psi(\tau)}(1+\tau)^\alpha|v(\tau)|^p}_2
 \norm{e^{\Psi(\tau)}u_\tau(\tau)}_2\notag\\
 &\leq
 \sqrt2\,\norm{e^{\Psi(\tau)}(1+\tau)^\alpha|v(\tau)|^p}_2
 \E_\Psi[u](\tau)^{\frac{1}{2}}.
\end{align}
It remains to estimate the exponentially weighted source. \Cref{lem:GN} with $q=2p$
and $\delta=\frac{1}{p}$ gives
\begin{align*}
 \norm{e^{\Psi(\tau)}|v(\tau)|^p}_2
 &=\norm{e^{\frac{\Psi(\tau)}{p}}v(\tau)}_{2p}^p \leq
 C(1+\tau)^{p(1-\vartheta(2p))}
 \norm{\nabla v(\tau)}_2^{p-1}
 \norm{e^{\Psi(\tau)}\nabla v(\tau)}_2.
\end{align*}
Using \eqref{eq:X-norm}, together with
$\vartheta(2p)=d\frac{p-1}{2p}$ and
$\alpha=\rho(p-1)-2$, we find
\begin{equation}
\label{eq:weighted-source}
 \norm{e^{\Psi(\tau)}(1+\tau)^\alpha|v(\tau)|^p}_2
 \leq
 C(1+\tau)^{-1-(d-\rho)(p-1)}
 \norm v_{X(T)}^p.
\end{equation}
Since $\rho<d$, the right-hand side is integrable over $[0,\infty)$.
To avoid dividing by an energy that may vanish, we apply
\eqref{eq:weighted-energy-differential} to
$(\E_\Psi[u]+\kappa)^{\frac{1}{2}}$ and let $\kappa$ tend to $0$. Using
\eqref{eq:weighted-source}, we obtain
\begin{align}
\label{eq:N-weighted-energy}
 \sup_{0\leq\tau\leq T}\E_\Psi[\N v](\tau)^{\frac{1}{2}}
 &\leq
 \E_\Psi[\N v](0)^{\frac{1}{2}}
 +C\int_0^T
 \norm{e^{\Psi(s)}(1+s)^\alpha|v(s)|^p}_2\dd s\notag\\ &\leq C\eps+C\norm v_{X(T)}^p.
\end{align}
It remains only to recover the exponentially  weighted $L^2$-norm of $\N v$.
\Cref{lem:GN} with $q=2$ and $\delta=1$ yields
\[
 (1+\tau)^{-1}\norm{e^{\Psi(\tau)}\N v(\tau)}_2
 \leq C\norm{e^{\Psi(\tau)}\nabla\N v(\tau)}_2.
\]
Consequently, \eqref{eq:N-weighted-energy} gives
\begin{align}
\label{eq:N-weighted-final}
 &\sup_{0\leq\tau\leq T}\left[
 (1+\tau)^{-1}\norm{e^{\Psi(\tau)}\N v(\tau)}_2
 +\norm{e^{\Psi(\tau)}
 (\nabla\N v,\partial_\tau\N v)(\tau)}_2
 \right]\notag\\
 &\qquad\leq C\eps+C\norm v_{X(T)}^p.
\end{align}
Combining \eqref{eq:N-unweighted-final} and
\eqref{eq:N-weighted-final}, we conclude that
\begin{equation}
\label{eq:N-bound}
 \norm{\N v}_{X(T)}
 \leq C\eps+C\norm v_{X(T)}^p,
\end{equation}
uniformly in $T$.

For completeness, let us establish the contraction estimate. Let
\(v,w\in X(T)\) and consider
\(
\mathcal{N}v-\mathcal{N}w.
\)
The key difference estimate is
\[
\bigl||v|^p-|w|^p\bigr|
\leq
C_p\bigl(|v|^{p-1}+|w|^{p-1}\bigr)|v-w|.
\]
Therefore, by H\"older's inequality, \Cref{lem:GN}, and the corresponding
estimates for the difference, the same arguments used above yield
\begin{align}\label{eq:N-difference}
    \left\|\mathcal{N}v-\mathcal{N}w\right\|_{X(T)}
\leq
C\left(
\left\|v\right\|_{X(T)}^{p-1}
+
\left\|w\right\|_{X(T)}^{p-1}
\right)
\left\|v-w\right\|_{X(T)}.
\end{align}
Choose $R=2C\eps$. If $\eps>0$ is sufficiently small,
\eqref{eq:N-bound} shows that $\N$ maps the closed ball
\[
 B_R(T):=\{v\in X(T):\norm v_{X(T)}\leq R\}
\]
into itself, while \eqref{eq:N-difference} shows that $\N$ is a
contraction on $B_R(T)$. The constants are independent of $T$.
Banach's fixed-point theorem gives a solution of \eqref{eq:damped-v} on every
$[0,T]$, and uniqueness makes these solutions compatible as $T$
varies. Thus $v$ is global and
\begin{equation}
\label{eq:v-decay}
 \norm{v(\tau)}_2\leq C\eps(1+\tau)^{-\frac{d}{2}},
 \qquad
 \norm{(\nabla v,v_\tau)(\tau)}_2
 \leq C\eps(1+\tau)^{-\frac{d}{2}-1}.
\end{equation}
The classical Gagliardo--Nirenberg inequality also gives
\begin{equation}
\label{eq:v-Lq}
 \norm{v(\tau)}_q
 \leq C_q\eps(1+\tau)^{-d(1-\frac{1}{q})},
 \qquad
 2\leq q\leq\frac{2d}{(d-2)_+}.
\end{equation}

Finally, $1+\tau=e^t$ and $\phi=e^{\rho t}v$ (see \eqref{eq:change}), while
\[
 \phi_t=e^{\rho t}\bigl((1+\tau)v_\tau+\rho v\bigr),
 \qquad
 \nabla\phi=e^{\rho t}\nabla v.
\]
The estimates \eqref{eq:v-decay} and \eqref{eq:v-Lq} imply the desired
estimates \eqref{eq:weighted-phi-energy} and
\eqref{eq:weighted-phi-Lq}. This completes the proof.
\end{proof}

%%%%%%%%%%
%%%%%%%%%%
%%%%%%%%%%

\subsection{Proof of  \Cref{thm:blowup}}

Finally, we turn to the proof of the blow-up result. The main idea is to reduce
the original partial differential equation to a nonlinear ordinary
differential inequality for the spatial average of the solution. More
precisely, the property of finite propagation speed of solutions allows us to test the weak formulation \eqref{eq:weak-formulation}
with a
function that is constant on the support of the solution. This yields a
scalar differential equation for the spatial average, while H\"older's
inequality provides a lower bound for its nonlinear forcing term. After
removing the unstable linear mode by a suitable exponential change of
unknown, the problem is reduced to the ODE criterion established below.

The next statement is an elementary quantitative version of the classical
ODE blow-up criterion introduced by Kato and widely used in the theory of
semilinear wave equations; see \cite{Kato1980,Sideris1984}. Quantitative
refinements designed to derive upper bounds for the lifespan can be found
in \cite{Takamura2015}. Related differential-inequality arguments for wave
equations with scale-invariant lower-order terms are employed in
\cite{NascimentoPalmieriReissig2017,TodorovaYordanov2001}. We include a complete proof because
the explicit dependence on the lower level \(\eta\) in
\eqref{eq:ode-lifespan} is needed below.

\begin{lemma}\label{lem:ode}
Let $a,b>0$, $p>1$, and $0<\eta\le1$. Suppose that
$y\in W_{\mathrm{loc}}^{2,1}([t_0,T))$ satisfies
\begin{equation}\label{eq:quantitative-ode}
 y(t)\ge c_0\eta,\qquad y'(t)\ge0,\qquad
 y''(t)+ay'(t)\ge by(t)^p
\end{equation}
for almost every $t\in[t_0,T)$, where $c_0>0$ is independent of $\eta$.
Then
\begin{equation}\label{eq:ode-lifespan}
 T-t_0\le C\eta^{-(p-1)},
\end{equation}
where $C$ depends only on $a,b,p$, and $c_0$. In particular, $T<\infty$.
\end{lemma}

\begin{proof}
Since $y\in W_{\mathrm{loc}}^{2,1}$, the function $y'$ is locally absolutely
continuous and all the integrations below are legitimate. Moreover, $y$ is
nondecreasing because $y'\ge0$ almost everywhere.

Set
\[
 M_0:=y(t_0),
 \qquad
 M_k:=2^kM_0 \quad (k\ge0).
\]
Due to \eqref{eq:quantitative-ode} it holds $M_0>0$.
We estimate the time required for $y$ to pass from the level $M_k$ to the
level $M_{k+1}=2M_k$. Whenever these levels are reached, define
\[
 t_k:=\inf \big\{t\in[t_0,T):y(t)\ge M_k\big\}.
\]
The estimates below also show that, whenever $t_k<T$ and the available
existence interval is longer than the resulting upper bound, the next level
$M_{k+1}$ is necessarily reached. Thus no existence of a hitting time is
being assumed in advance.
By the definition of $M_0$, the first hitting time is precisely the original
time $t_0$. By continuity and monotonicity,
\[
 y(t_k)=M_k,
 \qquad
 M_k\le y(t)<2M_k
 \quad\text{for }t_k\le t<t_{k+1}.
\]
Fix $k$ and write $t=t_k+\tau$. On the interval before the next level is
reached, \eqref{eq:quantitative-ode} gives
\[
 y''(t)+ay'(t)\ge bM_k^p.
\]
Multiplying by $e^{at}$ and integrating from $t_k$ to $t_k+\tau$, we obtain
\begin{align} \label{inequablowup}
 y'(t_k+\tau)
 &\ge e^{-a\tau}y'(t_k)
   +\frac baM_k^p\big(1-e^{-a\tau}\big) \ge \frac baM_k^p\big(1-e^{-a\tau}\big).
\end{align}
This estimate has two useful consequences, according to the size of $M_k$.

\medskip
\noindent
\emph{Step 1: levels below one.}
For $\tau\ge a^{-1}$,
\(
 1-e^{-a\tau}\ge1-e^{-1},
\)
and hence, by \eqref{inequablowup}
\begin{equation}\label{eq:parabolic-speed}
 y'(t_k+\tau)\ge c_1M_k^p,
 \qquad
 c_1:=\frac ba\big(1-e^{-1}\big).
\end{equation}
If $t_{k+1}-t_k\le a^{-1}$, there is nothing to prove. Otherwise, integrating
\eqref{eq:parabolic-speed} from $t_k+a^{-1}$ to $t_{k+1}$ gives
\begin{align*}
 c_1M_k^p(t_{k+1}-t_k-a^{-1})
 &\le y(t_{k+1})-y(t_k+a^{-1}) \le 2M_k-M_k=M_k.
\end{align*}
Consequently,
\begin{equation}\label{eq:low-doubling}
 t_{k+1}-t_k
 \le a^{-1}+c_1^{-1}M_k^{-(p-1)}.
\end{equation}

Let $K$ be the smallest nonnegative integer such that $M_K\ge1$. If
$M_0\ge1$, take $K=0$. For $0\le k<K$, we have $M_k<1$ and therefore,
$M_k^{-(p-1)}\ge1$. Thus the first term of the right-hand side in
\eqref{eq:low-doubling} can be absorbed into the second one:
\[
 t_{k+1}-t_k\le C_1M_k^{-(p-1)}.
\]
Summing the resulting geometric series gives
\begin{align}
 t_K-t_0
 &\le C_1\sum_{k=0}^{K-1}M_k^{-(p-1)}
 =C_1M_0^{-(p-1)}
   \sum_{k=0}^{K-1}2^{-k(p-1)}\notag\le \frac{C_1}{1-2^{-(p-1)}}M_0^{-(p-1)} \\
 &\le C_2c_0^{-(p-1)}\eta^{-(p-1)},
\label{eq:low-level-sum}
\end{align}
where the last inequality comes from
\eqref{eq:quantitative-ode}.

\medskip
\noindent
\emph{Step 2: levels above one.}
For $0\le z\le1$, the function $(1-e^{-z})/z$ is bounded from below by
$1-e^{-1}$. Hence, whenever $0\le\tau\le a^{-1}$,
\(
 1-e^{-a\tau}\ge (1-e^{-1})a\tau.
\)
It follows that by \eqref{inequablowup}
\begin{equation}\label{eq:wave-speed}
 y'(t_k+\tau)\ge c_2M_k^p\tau,
 \qquad
 c_2:=b(1-e^{-1}),
 \qquad 0\le\tau\le a^{-1}.
\end{equation}
Integrating once more, we find
\begin{equation}\label{eq:quadratic-growth}
 y(t_k+\tau)-M_k
 \ge \frac{c_2}{2}M_k^p\tau^2,
 \qquad 0\le\tau\le a^{-1}.
\end{equation}
Choose
\[
 A:=\max\Big\{1,\Big(\frac{2a^2}{c_2}\Big)^{\frac{1}{p-1}}\Big\}.
\]
If $M_k\ge A$, set
\[
 \tau_k:=\Big(\frac{2}{c_2}\Big)^{\frac{1}{2}}M_k^{-\frac{p-1}{2}}.
\]
The definition of $A$ ensures that $\tau_k\le a^{-1}$. Substituting
$\tau=\tau_k$ into \eqref{eq:quadratic-growth} gives
\[
 y(t_k+\tau_k)-M_k\ge M_k.
\]
Thus $y(t_k+\tau_k)\ge2M_k=M_{k+1}$, and the definition of $t_{k+1}$ implies
\[
 t_{k+1}\le t_k+\tau_k,
 \quad\text{and hence,}\quad
 t_{k+1}-t_k\le\tau_k
 =\Big(\frac{2}{c_2}\Big)^{\frac12}M_k^{-\frac{p-1}{2} }=C_3M_k^{-\frac{p-1}{2}}
 \quad(M_k\ge A).
\]
Only finitely many dyadic levels lie in $[1,A)$. Their number and, by
\eqref{eq:low-doubling}, the total time needed to cross them depend only on
$a,b,p$. Starting from the first level not smaller than $A$, we may sum the previous estimate:
\begin{align*}
 \sum_{M_k\ge A}(t_{k+1}-t_k)
 &\le C_3\sum_{M_k\ge A}M_k^{-\frac{p-1}{2}}\le C_4\sum_{j=0}^{\infty}2^{-j\frac{p-1}{2}}
 \le C_5.
\end{align*}
Combining this bound with \eqref{eq:low-level-sum}, and using
$\eta^{-(p-1)}\ge1$, we obtain
\begin{equation}\label{eq:total-time}
 \sum_{k=0}^{\infty}(t_{k+1}-t_k)
 \le C\eta^{-(p-1)}.
\end{equation}
Finally, we define
\[
d_k:=
\begin{cases}
 C_1M_k^{-(p-1)}, & M_k<1,\\[1mm]
 a^{-1}+c_1^{-1}M_k^{-(p-1)}, & 1\le M_k<A,\\[1mm]
 C_3M_k^{-\frac{p-1}{2}}, & M_k\ge A.
\end{cases}
\]
By the preceding estimates,
\[
t_{k+1}-t_k\le d_k.
\]
Splitting the sum according to the three ranges above and using the
same geometric-series estimates as before, we obtain
\[
\sum_{k=0}^{\infty}d_k
\le C\eta^{-(p-1)}.
\]
Consequently,
\[
t_k-t_0
\le \sum_{j=0}^{k-1}d_j
\le C\eta^{-(p-1)}
\qquad\text{for every }k\ge1.
\]
If $T-t_0>C\eta^{-(p-1)}$, then all the times $t_k$ would belong to
the compact interval
\(
 [t_0,t_0+C\eta^{-(p-1)}]\subset[t_0,T).
\)
However,
\(
 y(t_k)=M_k=2^kM_0\longrightarrow\infty,
\)
which contradicts the boundedness of the continuous function $y$ on
this compact interval. Therefore,
\[
 T-t_0\le C\eta^{-(p-1)},
\]
and, in particular, $T<\infty$.
\end{proof}

\begin{proof}[Proof of \Cref{thm:blowup}]
Fix a weak solution of \eqref{eq:phi} on $[0,T)$ in the sense of
\Cref{def:weak-solution}, where $T>0$. %\\

We prove an upper bound for this
arbitrary $T$. Then taking the supremum over all such intervals will give the
same bound for $T_\varepsilon$.

\medskip
\noindent
\emph{Step 1: equation for the spatial average.}
Fix $T_1<T$. By \eqref{eq:weak-support}, for almost every
$t\in[0,T_1]$ the function $\phi(t,\cdot)$ is supported in the fixed ball
$B_{R+e^{T_1}-1}$. Hence,
\[
 F(t):=\int_{\R^d}\phi(t,x)\dd x,
 \qquad
 H(t):=\int_{\R^d}|\phi(t,x)|^p\dd x
\]
belong to $L^1(0,T_1)$. Choose $\chi\in C_0^\infty(\R^d)$ such that
$\chi\equiv1$ on $B_{R+e^{T_1}-1}$ and insert
$\psi(t,x)=\zeta(t)\chi(x)$, with
$\zeta\in C_0^\infty([0,T_1))$, into
\eqref{eq:weak-formulation}. Since $\Delta\chi=0$ on the support of $\phi$,
we obtain
\begin{align*}
 \int_0^{T_1}F(t)
 \bigl(\zeta''(t)-n\zeta'(t)+m^2\zeta(t)\bigr)\dd t
 &=\int_0^{T_1}H(t)\zeta(t)\dd t+\varepsilon(F_1+nF_0)\zeta(0)
 -\varepsilon F_0\zeta'(0),
\end{align*}
where
\[
 F_0:=\int_{\R^d}f(x)\dd x,
 \qquad
 F_1:=\int_{\R^d}g(x)\dd x.
\]
Thus, in the sense of distributions on $(0,T_1)$,
\begin{equation}\label{eq:average-ode}
 F''+nF'+m^2F=H,
 \qquad
 F(0)=\varepsilon F_0,
 \qquad
 F'(0)=\varepsilon F_1.
\end{equation}
Since $F,H\in L^1(0,T_1)$, the standard one-dimensional regularity for this
constant-coefficient equation yields $F\in W^{2,1}(0,T_1)$. Because $T_1<T$
was arbitrary, we may conclude
\[
 F\in W_{\mathrm{loc}}^{2,1}([0,T)).
\]
In particular, the identities below hold almost everywhere and their
integrated forms hold at every time.

\noindent The propagation property also gives
\[
 |\operatorname{supp}\phi(t,\cdot)|
 \le |B_{R+e^t-1}|\le C_R e^{dt}.
\]
Therefore, H\"older's inequality implies
\begin{equation}\label{eq:H-lower}
 H(t)
 \ge |\operatorname{supp}\phi(t,\cdot)|^{1-p}|F(t)|^p
 \ge C_R e^{-d(p-1)t}|F(t)|^p
\end{equation}
for almost every $t\in[0,T)$.

\medskip
\noindent
\emph{Step 2: removal of the unstable linear mode.}
Since \(m^2<0\), the parameters introduced in
\eqref{parameternu} and \eqref{def:rho} take the form
\[
 \nu:=\sqrt{\frac{n^2}{4}-m^2},
 \qquad
 \rho:=\nu-\frac n2.
\]
The hypothesis \eqref{restonmublowup} is equivalent to $\rho \ge d$. The two roots
of $\lambda^2+n\lambda+m^2=0$ are $\rho$ and $-(n+\rho)$. Variation of
constants in \eqref{eq:average-ode} gives
\begin{equation}\label{eq:F-representation}
 F(t)=A_\varepsilon e^{\rho t}
      +B_\varepsilon e^{-(n+\rho)t}
      +\frac1{2\nu}\int_0^t
       \bigl(e^{\rho(t-s)}-e^{-(n+\rho)(t-s)}\bigr)H(s)\dd s,
\end{equation}
where
\begin{equation}\label{eq:A-B}
 A_\varepsilon
 =\frac{\varepsilon}{2\nu}\bigl(F_1+(n+\rho)F_0\bigr)>0,
 \qquad
 B_\varepsilon
 =\frac{\varepsilon}{2\nu}(\rho F_0-F_1).
\end{equation}
The integral kernel in \eqref{eq:F-representation} is nonnegative. Moreover,
$A_\varepsilon=\varepsilon A_1$ and
$B_\varepsilon=\varepsilon B_1$, where $A_1>0$ and $B_1$ are independent of
$\varepsilon$. We can therefore choose $t_0\ge0$, independently of
$\varepsilon$, such that
\[
 |B_1|e^{-(n+2\rho)t_0}\le\frac{A_1}{2}.
\]
It follows from \eqref{eq:F-representation} that
\begin{equation}\label{eq:F-positive}
 F(t)\ge\frac{A_\varepsilon}{2}e^{\rho t}>0,
 \qquad t\ge t_0.
\end{equation}
If $T\le t_0$, both estimates in \eqref{eq:lifespan} follow after
enlarging $C$. We may therefore assume below that $T>t_0$.

Define
\(
 Y(t):=e^{-\rho t}F(t).
\)
Since $\rho$ is a characteristic root, \eqref{eq:average-ode} and
\eqref{eq:H-lower} yield
\begin{align}\label{eq:Y-main}
 Y''+2\nu Y'
 =e^{-\rho t}H(t) \ge C_R e^{(p-1)(\rho-d)t}|Y(t)|^p.
\end{align}
Furthermore, \eqref{eq:F-positive} gives
\begin{equation}\label{eq:Y-lower}
 Y(t)\ge\frac{A_1}{2}\varepsilon=:c_0\varepsilon,
 \qquad t\ge t_0.
\end{equation}
In particular, \eqref{eq:Y-main} implies
\begin{equation}\label{eq:Y-forcing}
 Y''+2\nu Y'
 \ge c\varepsilon^p e^{(p-1)(\rho-d)t},
 \qquad t\ge t_0.
\end{equation}
The representation formula also gives $|Y'(t_0)|\le C\varepsilon$.

\medskip
\noindent
\emph{Step 3: the case $\rho>d$.}
Let $\sigma:=\rho-d>0$. Multiplication of
\eqref{eq:Y-forcing} by $e^{2\nu t}$ and integration from $t_0$ to $t$ give
\begin{equation}\label{eq:Yprime-supercritical}
 Y'(t)
 \ge-C\varepsilon e^{-2\nu(t-t_0)}
 +c\varepsilon^p e^{(p-1)\sigma t}
 \qquad(t\ge t_0+1),
\end{equation}
after changing the constants.

Set
\(
t_\varepsilon^*:=\frac1\sigma\log(\varepsilon^{-1}).
\)
Provided that the interval of existence reaches this time, the positive term on the right-hand side of
\eqref{eq:Yprime-supercritical} is of order $\varepsilon$, whereas the
negative term is $O(\varepsilon^{1+\frac{2\nu}{\sigma}})$. Hence, after adding a
constant independent of $\varepsilon$, we obtain a time
\begin{equation}\label{eq:teps-supercritical}
 s_\varepsilon\le
 \frac1\sigma\log(\varepsilon^{-1})+C
\end{equation}
such that
\[
 Y(s_\varepsilon)\ge c_0\varepsilon,
 \qquad
 Y'(s_\varepsilon)\ge c_1\varepsilon.
\]
If $T\le s_\varepsilon$, estimate \eqref{eq:lifespan} is already true.
Assume therefore that $T>s_\varepsilon$ and define
\[
 Z(\tau):=\varepsilon^{-1}Y(s_\varepsilon+\tau),
 \qquad 0\le\tau<T-s_\varepsilon.
\]
Using \eqref{eq:Y-main}, \eqref{eq:teps-supercritical}, and
$e^{(p-1)\sigma s_\varepsilon}\ge c\varepsilon^{-(p-1)}$, we obtain
\begin{equation}\label{eq:Z-ode}
 Z''+2\nu Z'
 \ge b e^{(p-1)\sigma\tau}|Z|^p
 \ge b|Z|^p,
\end{equation}
where $b>0$ is independent of $\varepsilon$.

We now verify the monotonicity required by \Cref{lem:ode}. Initially
$Z(0)\ge c_0$ and $Z'(0)\ge c_1>0$. As long as $Z$ is positive,
\eqref{eq:Z-ode} gives
\[
 \frac{\dd}{\dd\tau}\bigl(e^{2\nu\tau}Z'(\tau)\bigr)\ge0.
\]
Thus $Z'(\tau)>0$, so $Z(\tau)\ge Z(0)>0$. By continuity, this argument
extends to the whole interval $[0,T-s_\varepsilon)$. Consequently,
\[
 Z(\tau)\ge c_0,
 \qquad Z'(\tau)\ge0,
 \qquad Z''+2\nu Z'\ge bZ^p.
\]
\Cref{lem:ode}, applied with $\eta=1$, yields
\(
 T-s_\varepsilon\le C.
\)
Together with \eqref{eq:teps-supercritical}, this proves
\eqref{eq:lifespan}.

\medskip
\noindent
\emph{Step 4: the endpoint $\rho=d$.}
In this case \eqref{eq:Y-forcing} becomes
\[
 Y''+2\nu Y'\ge c\varepsilon^p,
 \qquad t\ge t_0.
\]
After multiplication by $e^{2\nu t}$ and integration, we obtain
\begin{equation}\label{eq:Yprime-endpoint}
 Y'(t)
 \ge-C\varepsilon e^{-2\nu(t-t_0)}
 +c\varepsilon^p\bigl(1-e^{-2\nu(t-t_0)}\bigr).
\end{equation}
Provided that the interval of existence reaches the time on the right-hand
side below, it follows that there exists
\begin{equation}\label{eq:teps-endpoint}
 s_\varepsilon\le C\bigl(1+\log(\varepsilon^{-1})\bigr)
\end{equation}
such that $Y'(s_\varepsilon)\ge0$. If $T\le s_\varepsilon$, the desired
bound follows at once from
$1+\log(\varepsilon^{-1})\le C\varepsilon^{-(p-1)}$. Otherwise,
\eqref{eq:Y-main} and the integrating-factor argument show that
$Y'(t)\ge0$ for every $t\in[s_\varepsilon,T)$. Combining this with
\eqref{eq:Y-lower}, we have on this interval
\[
 Y(t)\ge c_0\varepsilon,
 \qquad
 Y'(t)\ge0,
 \qquad
 Y''(t)+2\nu Y'(t)\ge C_RY(t)^p.
\]
\Cref{lem:ode}, now with $\eta=\varepsilon$, gives
\(
 T-s_\varepsilon\le C\varepsilon^{-(p-1)}.
\)
Using \eqref{eq:teps-endpoint} once more, we conclude \eqref{eq:lifespan}.

We have proved \eqref{eq:lifespan} for every interval on which a weak
solution in the sense of \Cref{def:weak-solution} exists. Taking the
supremum over all such intervals proves the stated estimate for
$T_\varepsilon$ and completes the proof.
\end{proof}

\begin{remark}\label{concludingremark}
In a forthcoming paper we shall study as a generalization of \eqref{eq:phi}
the model
\begin{equation} \label{ourmodel} \begin{cases} \phi_{tt} - a(t)^2 \Delta \phi + n\frac{a(t)}{A(t)}\phi_t + m^2\Big(\frac{a(t)}{A(t)}\Big)^2 \phi=|\phi|^p, &(t,x)\in[0,\infty)\times\R^d,\\ \phi(0,x)=f(x),\qquad \phi_t(0,x)=g(x),&x\in\R^d,
 \end{cases}
\end{equation}
Here $A(t)=1+\int_0^t a(\tau)d\tau$, $n$ is a positive constant and $m^2$ is a real number.  This model \eqref{ourmodel} contains as a special case the anti-de Sitter model $(a(t)=e^{t})$. We are interested to study the influence of $m^2$ on global (in time) of small data well-posedness or blow up results, respectively. These considerations would continue those one of the very recent paper \cite{MarquesReissig}.
\end{remark}

% \section{Concluding remarks}

% The threshold has the same interpretation on both sides.  Above it, the
% Fourier-volume gain compensates for the unstable zero mode and makes the
% Duhamel term integrable.  At and below it, the zero mode grows at least as
% fast as the radius-induced volume scale and the spatial average obeys an
% autonomous superlinear inequality after conjugation.

% The positivity assumption is used only in the blow-up theorem.  It remains
% natural to ask whether sign-changing data at the endpoint admit a stable
% manifold, and whether the energy-subcritical restriction can be removed by
% working at higher Sobolev regularity.  These questions do not change the
% mass balance identified here.


\begin{thebibliography}{99}
\bibitem{Bachelot2011}
A.Bachelot,
\newblock The Klein--Gordon equation in the anti-de Sitter cosmology,
\newblock \emph{Journal de Math\'ematiques Pures et Appliqu\'ees} \textbf{96} (2011), 527--554.

\bibitem{Baskin2013}
D.Baskin,
\newblock Strichartz estimates on asymptotically de Sitter spaces,
\newblock \emph{Annales Henri Poincar\'e} \textbf{14} (2013), 221--252.

\bibitem{BreitenlohnerFreedman1982a}
P.Breitenlohner and D.Z.Freedman,
\newblock Positive energy in anti-de Sitter backgrounds and gauged extended supergravity,
\newblock \emph{Physics Letters B} \textbf{115} (1982), 197--201.

\bibitem{BreitenlohnerFreedman1982b}
P.Breitenlohner and D.Z.Freedman,
\newblock Stability in gauged extended supergravity,
\newblock \emph{Annals of Physics} \textbf{144} (1982), 249--281.

\bibitem{BuiReissig2014}
T.A.Bui and M.Reissig,
\newblock The interplay between time-dependent speed of propagation and dissipation in wave models,
\newblock in \emph{Fourier Analysis}, Trends in Mathematics, Birkh\"auser, 2014, 9--45.

\bibitem{DLMF}
F.W.J.Olver, A.B.Olde Daalhuis, D.W.Lozier, B.I.Schneider,
R.F.Boisvert, C.W.Clark, B.R.Miller, B.V.Saunders,
H.S.Cohl and M.A.McClain (eds.),
\newblock \emph{NIST Digital Library of Mathematical Functions},
\newblock National Institute of Standards and Technology,
\url{https://dlmf.nist.gov/}, Release 1.2.7, 2026.

\bibitem{EbertNascimento2018}
M.R.Ebert and W.Nunes do Nascimento,
\newblock A classification for wave models with time-dependent potential and speed of propagation,
\newblock \emph{Advances in Differential Equations} \textbf{23} (2018), no.11--12, 793--912.

\bibitem{EbertReissig2018}
M.R.Ebert and M.Reissig,
\newblock Regularity theory and global existence of small data solutions to semilinear de Sitter models with power non-linearity,
\newblock \emph{Nonlinear Analysis: Real World Applications} \textbf{40} (2018), 14--54.

\bibitem{Galstian2003}
A.Galstian,
\newblock $L^p$--$L^q$ decay estimates for wave equations with exponentially growing speed of propagation,
\newblock \emph{Applicable Analysis} \textbf{82} (2003), no.3, 197--214.

\bibitem{GalstianYagdjian2020}
A.Galstian and K.Yagdjian,
\newblock Finite lifespan of solutions of the semilinear wave equation in the Einstein--de Sitter spacetime,
\newblock \emph{Reviews in Mathematical Physics} \textbf{32} (2020), no.7, 2050018.

\bibitem{HintzVasy2015}
P.Hintz and A.Vasy,
\newblock Semilinear wave equations on asymptotically de Sitter,
Kerr--de Sitter and Minkowski spacetimes,
\newblock \emph{Analysis \& PDE} \textbf{8} (2015), no.8, 1807--1890.

\bibitem{Holzegel2012}
G.Holzegel,
\newblock Well-posedness for the massive wave equation on asymptotically anti-de Sitter spacetimes,
\newblock \emph{Journal of Hyperbolic Differential Equations} \textbf{9} (2012), no.2, 239--261.

\bibitem{HolzegelSmulevici2012}
G.Holzegel and J.Smulevici,
\newblock Self-gravitating Klein--Gordon fields in asymptotically anti-de Sitter spacetimes,
\newblock \emph{Annales Henri Poincar\'e} \textbf{13} (2012), no.4, 991--1038.

\bibitem{HolzegelWarnick2014}
G.Holzegel and C.M.Warnick,
\newblock Boundedness and growth for the massive wave equation on asymptotically anti-de Sitter black holes,
\newblock \emph{Journal of Functional Analysis} \textbf{266} (2014), no.4, 2436--2485.

\bibitem{IkehataTanizawa}
R.Ikehata and K.Tanizawa,
\newblock Global existence of solutions for semilinear dampedwave equations in $R^N$ with noncompactly supported initial data,
\newblock \emph{Nonlinear Anal} \textbf{61} (2005), 1189--1208.

\bibitem{IshibashiWald2004}
A.Ishibashi and R.M.Wald,
\newblock Dynamics in non-globally-hyperbolic static spacetimes. III. Anti-de Sitter spacetime,
\newblock \emph{Classical and Quantum Gravity} \textbf{21} (2004), 2981--3014.

\bibitem{MarquesReissig}
J.Marques and M.Reissig,
\newblock Global (in time) well-posedness of small data Sobolev solutions to
semilinear hyperbolic models with power nonlinearity and integrable decaying speed of propagation in the critical case,
\newblock accepted for publication in \emph{J. Differential Equations}, 41A4.

\bibitem{Nakamura2014}
M.Nakamura,
\newblock The Cauchy problem for semi-linear Klein--Gordon equations in de Sitter spacetime,
\newblock \emph{Journal of Mathematical Analysis and Applications} \textbf{410} (2014), no.1, 445--454.

\bibitem{Nakamura2015}
M.Nakamura,
\newblock Energy solutions for nonlinear Klein--Gordon equations in de Sitter spacetime,
\newblock in \emph{Current Trends in Analysis and Its Applications},
Trends in Mathematics, Birkh\"auser, Cham, 2015, 203--208.

\bibitem{Nakamura2021}
M.Nakamura,
\newblock The Cauchy problem for the Klein--Gordon equation under the quartic potential in the de Sitter spacetime,
\newblock \emph{Journal of Mathematical Physics} \textbf{62} (2021), no.12, 121509.

\bibitem{NascimentoPalmieriReissig2017}
W.Nunes do Nascimento, A.Palmieri and M.Reissig,
\newblock Semi-linear wave models with power non-linearity and
scale-invariant time-dependent mass and dissipation,
\newblock \emph{Mathematische Nachrichten} \textbf{290} (2017), 1779--1805.

\bibitem{Kato1980}
T.Kato,
\newblock Blow-up of solutions of some nonlinear hyperbolic equations,
\newblock \emph{Communications on Pure and Applied Mathematics}
\textbf{33} (1980), 501--505.

\bibitem{Palmieri2018}
A.Palmieri,
\newblock Global existence of solutions for semi-linear wave equation with
scale-invariant damping and mass in exponentially weighted spaces,
\newblock \emph{Journal of Mathematical Analysis and Applications}
\textbf{461} (2018), 1215--1240.

\bibitem{PalmieriTakamura2022}
A.Palmieri and H.Takamura,
\newblock A note on blow-up results for semilinear wave equations in de Sitter and anti-de Sitter spacetimes,
\newblock \emph{Journal of Mathematical Analysis and Applications} \textbf{514} (2022), no.1, 126266.

\bibitem{PalmieriTakamura2022critical}
A.Palmieri and H.Takamura,
\newblock On a semilinear wave equation in anti-de Sitter spacetime: the critical case,
\newblock \emph{Journal of Mathematical Physics} \textbf{63} (2022), no.11, 111505.

\bibitem{Sideris1984}
T.C.Sideris,
\newblock Nonexistence of global solutions to semilinear wave equations
in high dimensions,
\newblock \emph{Journal of Differential Equations}
\textbf{52} (1984), 378--406.

\bibitem{Takamura2015}
H.Takamura,
\newblock Improved Kato's lemma on ordinary differential inequality and
its application to semilinear wave equations,
\newblock \emph{Nonlinear Analysis}
\textbf{125} (2015), 227--240.

\bibitem{TodorovaYordanov2001}
G.Todorova and B.Yordanov,
\newblock Critical exponent for a nonlinear wave equation with damping,
\newblock \emph{Journal of Differential Equations}
\textbf{174} (2001), 464--489.

\bibitem{TsutayaWakasugi2020}
K.Tsutaya and Y.Wakasugi,
\newblock Blow up of solutions of semilinear wave equations in Friedmann--Lema\^itre--Robertson--Walker spacetime,
\newblock \emph{Journal of Mathematical Physics} \textbf{61} (2020), no.9, 091503.

\bibitem{TsutayaWakasugi2022}
K.Tsutaya and Y.Wakasugi,
\newblock Blow up of solutions of semilinear wave equations related to nonlinear waves in de Sitter spacetime,
\newblock \emph{Partial Differential Equations and Applications} \textbf{3} (2022), Paper No.6.

\bibitem{Vasy2012}
A.Vasy,
\newblock The wave equation on asymptotically anti-de Sitter spaces,
\newblock \emph{Analysis \& PDE} \textbf{5} (2012), no.1, 81--144.

\bibitem{Warnick2013}
C.M.Warnick,
\newblock The massive wave equation in asymptotically AdS spacetimes,
\newblock \emph{Communications in Mathematical Physics} \textbf{321} (2013), 85--111.

\bibitem{WatsonBessel}
G.N.Watson,
\newblock \emph{A Treatise on the Theory of Bessel Functions}, 2nd ed.,
\newblock Cambridge University Press, Cambridge, 1944.

\bibitem{Yagdjian2009}
K.Yagdjian,
\newblock The semilinear Klein--Gordon equation in de Sitter spacetime,
\newblock \emph{Discrete and Continuous Dynamical Systems Series S} \textbf{2} (2009), no.3, 679--696.

\bibitem{Yagdjian2010}
K.Yagdjian,
\newblock Fundamental solutions for hyperbolic operators with variable coefficients,
\newblock \emph{Rendiconti dell'Istituto di Matematica dell'Universit\`a di Trieste} \textbf{42} (2010), suppl., 221--243.

\bibitem{Yagdjian2012}
K.Yagdjian,
\newblock Global existence of the scalar field in de Sitter spacetime,
\newblock \emph{Journal of Mathematical Analysis and Applications} \textbf{396} (2012), no.1, 323--344.

\bibitem{Yagdjian2013}
K.Yagdjian,
\newblock Huygens' principle for the Klein--Gordon equation in the de Sitter spacetime,
\newblock \emph{Journal of Mathematical Physics} \textbf{54} (2013), no.9, 091503.

\bibitem{Yagdjian2015}
K.Yagdjian,
\newblock Integral transform approach to solving Klein--Gordon equation with variable coefficients,
\newblock \emph{Mathematische Nachrichten} \textbf{288} (2015), no.17--18, 2129--2152.

\bibitem{Yagdjian2019}
K.Yagdjian,
\newblock Global existence of the self-interacting scalar field in the de Sitter universe,
\newblock \emph{Journal of Mathematical Physics} \textbf{60} (2019), no.5, 051503.

\bibitem{YagdjianGalstian2008}
K.Yagdjian and A.Galstian,
\newblock Fundamental solutions of the wave equation in Robertson--Walker spaces,
\newblock \emph{Journal of Mathematical Analysis and Applications} \textbf{346} (2008), no.2, 501--520.

\bibitem{YagdjianGalstian2009dS}
K.Yagdjian and A.Galstian,
\newblock Fundamental solutions for the Klein--Gordon equation in de Sitter spacetime,
\newblock \emph{Communications in Mathematical Physics} \textbf{285} (2009), 293--344.

\bibitem{YagdjianGalstian2009AdS}
K.Yagdjian and A.Galstian,
\newblock The Klein--Gordon equation in anti-de Sitter spacetime,
\newblock \emph{Rendiconti del Seminario Matematico della Universit\`a e del Politecnico di Torino}
\textbf{67} (2009), no.2, 271--292.
\end{thebibliography}
\end{document}